\documentclass {article}

\usepackage{PRIMEarxiv}

\newcommand{\ra}[1]{\renewcommand{\arraystretch}{#1}}

\usepackage[
singlelinecheck=false 
]{caption}
\usepackage{algorithm}
\usepackage[]{algpseudocode}
\algdef{SE}[SUBALG]{Indent}{EndIndent}{}{\algorithmicend\ }%
\algtext*{Indent}
\algtext*{EndIndent}
\usepackage{amsmath} 
\usepackage{amssymb}
\allowdisplaybreaks 
\usepackage{breakurl}
\usepackage[breaklinks, hidelinks]{hyperref} 
\usepackage{booktabs}
\usepackage{color}
\usepackage{color, colortbl}
\definecolor{Gray}{gray}{0.9}

\usepackage{enumitem}
\usepackage{eqparbox}
\usepackage{eurosym}
\usepackage{gensymb}

\usepackage{graphicx}
\usepackage[labelsep=space]{caption}
\usepackage{longtable}
\usepackage{mwe}
\usepackage{multirow}
\usepackage{natbib}%

\usepackage[figuresright]{rotating}
\usepackage{setspace}
\usepackage{siunitx}

\usepackage{subcaption}
\usepackage{soul}
\usepackage{tikz}
\usetikzlibrary{arrows,shapes,snakes,automata,backgrounds,petri,fit,positioning, chains}
\usepackage{url}
\usepackage{xfrac}

\newcommand\correspondingauthor{\thanks{Corresponding author}}

\title{Enhancing Pharmaceutical Cold Supply Chain: Integrating Medication Synchronization and Diverse Delivery Modes
}

\author{
Elise Potters\\
Department of Industrial Engineering \\ and Business Information Systems \\
University of Twente, Enschede, Netherlands \\ 
  \And
Behzad Mosalla Nezhad \\
   Department of Industrial Engineering \\ and Business Information Systems \\
  University of Twente, Enschede, Netherlands \\ 
   \And
  Viktor Huiskes\\
  Maartensapotheek, Sint Maartenskliniek \\
  Nijmegen,  Netherlands \\ 
  \And
    Erwin Hans\\
     Department of Industrial Engineering \\ and Business Information Systems \\
  University of Twente, Enschede, Netherlands \\
  \And
   Amin Asadi\correspondingauthor \\
  University of Twente, Enschede, Netherlands \\
   Department of Industrial Engineering \\ and Business Information Systems \\
  \texttt{amin.asadi@utwente.nl} 
}

\begin{document}
\maketitle
\begin{abstract}
The significance of last-mile logistics in the healthcare supply chain is growing steadily, especially in pharmacies where the growing prevalence of medication delivery to patients' homes is remarkable. This paper proposes a novel mathematical model for the last-mile logistics of the pharmaceutical supply chain and optimizes a pharmacy's logistical financial outcome while considering medication synchronization, different delivery modes, and temperature requirements of medicines. We propose a mathematical formulation of the problem using Mixed Integer Linear Programming (MILP) evolved from the actual problem of an outpatient pharmacy of a Dutch hospital. We create a case study by gathering, preparing, processing, and analyzing the associated data. We find the optimal solution, using Python MIP package and the Gurobi solver, which indicates the number of order batches, the composition of these batches, and the number of staff related to the preparation of the order batches. Our results show that our optimal solution increases the pharmacy's logistical financial outcome by 34 percent. Moreover, we propose other model variations and perform extensive scenario analysis to provide managerial insights applicable to other pharmacies and distributors in the last step of cold supply chains. Based on our scenario analysis, we conclude that improving medication synchronization can significantly enhance the pharmacy's logistical financial outcome.
\end{abstract}

\keywords{Supply Chain Management \and Mixed Integer Linear Programming \and
Healthcare Delivery \and Medication Synchronization}

\section{Introduction and Background}\label{sec:Int}

The global supply chain management (SCM) market in healthcare was valued at \$2.33 billion in 2021 and is expected to expand at a compound annual growth rate of 9.2\% from 2022 to 2030 \citep{Grandviewresearch2020Healthcare2030}. Within the healthcare supply chain, the importance of last-mile logistics in the pharmaceutical sector is rapidly growing \citep{Sun2022COVIDSystem, Dcruz2022TheAhead} as pharmacies tend to provide medication delivery to the patients' homes. While home delivery of medications can enhance patients' accessibility to pharmacy care and reduce their visits to the pharmacy, it presents complexities that require special consideration in the supply chain \citep{Srivastava2022TemperatureProducts}. Examples of such complications include implementing additional safety measures related to specific temperature and storage requirements for certain pharmaceuticals during transportation and delivery \citep{Srivastava2022TemperatureProducts}. As a result, a growing number of researchers, including ourselves, are focusing on managing the distribution operations of these cold supply chains. One more complex yet advantageous element that can be integrated into pharmaceutical supply chains (PSCs) is the synchronization or consolidation of medications. Optimal medication synchronization can reduce workload and costs \citep{Melendez2020ReducingCenters} and increase medication adherence \citep{Ross2013SyncSwim}. In this research, We aim to explore the implementation of medication synchronization while considering different delivery modes and adhering to the temperature requirements of medicines in transportation and study their impact on the logistical financial outcome of the pharmacy.

The pharmaceutical supply chain network (PSCN) is a complex structure that involves the manufacturing and distributing of medications to patients and organizations at the right time, quantity, and cost \citep{Abdulghani2019PharmaceuticalDrugs}. Managing PSCNs requires addressing challenges in logistics, inventory management, and warehouse management \citep{Sinnei2023EvaluationKenya}. In recent years, cost optimization  
has become a cross-sectional theme in the supply chains \citep{FitchSolutions2022InflationEurope}. Especially with the emergence of the COVID-19 outbreak, well-organized cold supply chain networks are perceived to be crucial more than ever given the high global demand and varying temperature requirements for such vaccines and medicines \citep{Balfour2021NearlyCOVID-19}. 
Failure to maintain the recommended temperature range at one or multiple stages of the supply chain can lead to decreased effectiveness or even harm for various pharmaceuticals, which is not limited to the COVID-19 vaccine case \citep{Nyirimanzi2023AssessmentRwanda, Turan2022InvestigatingIndustry}. This can result in waste of supplies, increased medical costs, and increased disease prevalence, as evidenced by previous research \citep{Malik2022StakeholdersStudy} highlighting the necessity of effective cold supply chain management in pharmacy logistics \citep{Bishara2006ColdChain}. Therefore, this research aims to provide a framework for coping with these challenges. Examples of research toward studying SCM (particularly cold SCM) from other close perspectives can be found in \citep{Lagana2022MultipleResearch, Dixit2019AResearch}. \cite{Kartoglu2014ToolsChain} presents several tools and approaches to ensure the quality of vaccines throughout the cold chain wherein vaccines require adequate temperature storage, as they are either sensitive to heat or cold. Similar to our work, this research focuses on cold SCM with consideration of cost factors and different delivery systems. However, they have different approaches and aim to improve the cold PSCN by avoiding the negative consequences of faulty management.


Last-mile logistics, the final stage of the supply chain where a shipment is carried to its ultimate destination, plays an indispensable role in delivering medications safely and timely to patients that rigorously rely on different delivery modes \citep{Giuffrida2022OptimizationReview, Bhatnagar2018CorrespondingRecommendations}. Here, we also explore various delivery modes for last-mile logistics within the PSCN. As aforementioned, the PSCN involves coordinating all activities related to the sourcing, procurement, manufacturing, and delivery of products and services in the pharmaceutical setting. There has been a growing interest in innovative delivery modes from the pharmacy to the patient, including pick-up services, bike couriers, electric vehicles, and even drones to demonstrate potential benefits like reduced carbon emissions and improved demand satisfaction \citep{Lamiscarre2022AssessingLogistics, Asadi2022ADeliveries}. For instance, \cite{Handoko2014AnCentre} proposed a profit-maximizing model for determining which demands to serve using an Urban Consolidation Centre (UCC). Unlike this work, our research studies the financial impact of various delivery strategies, focusing on the various last-mile delivery modes within the PSCN. The use of crowd shipping for delivery is discussed in \cite{Bajec2022ADeliveries}, highlighting its benefits (low prices, reduced delivery times, and environmental impact) and risks (lower quality, reliability, and inefficiency), which can be mitigated by proper management of coalition members and optimal collaboration. Another relevant study is \cite{Azizi2022ReliableFactor}, which proposed and optimized a hub-and-spoke system for transportation between multiple locations. This system utilizes a hub to link several locations, possibly other hubs, in a distribution network. The study focuses on determining the optimal location of the hubs and demand allocation to all hubs. Moreover, \cite{Lamiscarre2022AssessingLogistics} evaluated the potential of drone-based last-mile logistics as a modern delivery mode that could supplement or replace traditional delivery options. The study explores a hybrid approach that combines both ground vehicles and Unmanned Aerial Vehicles (UAVs) and assesses the benefits of this approach in terms of time, pollution, and cost. This research is relevant to our investigation of multiple delivery modes. Besides, we anticipate that the use of UAVs will become more prevalent soon, and \cite{Lamiscarre2022AssessingLogistics} provides valuable insights into the potential use of UAVs. Our study does not include UAVs as this is not a feasible solution for the case study at the MA. However, we consider e-bike delivery, which is a modern, low-cost, and sustainable mode.  

Besides delivery modes, medication synchronization or medication consolidation is a crucial aspect of optimizing the delivery of pharmaceuticals, which refers to the action of filling patients’ medication prescriptions simultaneously \citep{Nguyen2017TheReview,Pham2023PilotProgram}. It undoubtedly leads to aligning patients' medication refill dates and reducing the need for frequent pharmacy visits \citep{Melendez2020ReducingCenters}. This strategy has the potential to enhance medication adherence, decrease environmental impacts and costs, and improve the patient experience \citep{Ross2013SyncSwim}. \cite{Ross2013SyncSwim} highlights the importance of refill consolidation, which is the synchronization of refill dates for patients who use multiple medications. They surveyed 50 patients and reported that poor refill consolidation is linked to reduced medication adherence, and it could be harmful to patients. Furthermore, \cite{Hughes2022AnAnalysis} explored practical solutions for the barriers and facilitators associated with Med Sync (a medication synchronization service) in community pharmacies, such as engaging staff, organizing programs effectively, and collaborating with providers. However, the financial outcome of medication synchronization in the pharmaceutical sector needs further study. 


In the pharmaceutical industry, where margins are typically thin \citep{Garattini2008PricesCountries}, positive net financial outcomes are essential for the survival of these organizations. Thus, optimizing the cost of cold PSCN is of great significance. We note that considering factors such as temperature requirements, different delivery modes, and medication synchronization can add extra layers of complexity to the optimization problem. However, incorporating these factors into a proper mathematical framework 
can benefit pharmacies and other healthcare institutions. For instance, \cite{li2023two} presents a cold PSCN considering essential factors such as drug demand and safety risks. \cite{abdolazimi2023development} study mathematical modeling to control a non-cold PSCN in several directions, including costs, environmental impacts, lead time, and healthcare worker infection risk. Besides, they consider COVID-19 infection rates among medical staff and tackle the economic and environmental consequences. \cite{Chowdhury2022ModelingSystem} develop and optimize a vaccine supply chain, considering three objectives: minimizing cost, ensuring environmental and social sustainability, and maximizing job opportunities. Furthermore, \cite{Ala2023AnApproach} utilize the queuing theory and MILP framework to optimize waiting time and patient arrivals at a minimum supply chain cost. A study by \cite{Aghazadeh2018RobustChain} proposes a multi-objective MIP model for an organ transplant, aiming to reduce working center costs and optimize allocations of units and organs.

We note there are other optimization models, including non-linear programming \citep{Starita2022ImprovingInvestment, Juned2022DesigningNetwork}, simulation models \citep{Dui2023CascadingCOVID-19, Borges2020Simulation-basedHospital}, and Markov decision processes \citep{Abbaspour2021AUncertainty, Shakya2022ADemand}, which are widely used in healthcare SCM. Moreover, other considerations, such as sustainability have been embedded in PSCN problems, with studies by \cite{ahmad2022multi}, \cite{nasrollahi2021mathematical}, and \cite{halim2019decision} underscoring the significance of involving social and environmental concepts. An example of metaheuristic algorithms to optimize PSCN reliability can be found in \citep{goodarzian2020multi}. There are other relevant research discussions on the optimization models for PSCN \cite{zhu2018design, sousa2005global}, product launch operations \cite{hansen2015planning}, integration of procurement, production, and distribution decisions \cite{susarla2012integrated}. In contrast to the previous studies, the present work stands out by incorporating novel concepts and features within the PSCN model, such as medication synchronization, clustering each patient type, focusing on delivery mode, and staffing requirements.

The summary and comparison of the most current works can be found in Table \ref{PSCNLR}. Based on our literature review that shapes our knowledge, despite the abundance of studies on the PSCN, no research embedded medication synchronization and the incorporation of multiple cooling modes within PSCN settings. Hence, this work attempts to fill the existing gap and open a new avenue for future research in this area. Therefore, this study proposes the Pharmaceutical Cold Supply Chain Management Considering Medication Synchronization and Different Delivery Modes (PCSCM-MSDM) problem by exploring the impact of delivery methods, medication synchronization, and cooling requirements on distribution costs and revenues. The model optimizes the number of order batches and the composition of these batches, as well as the staffing related to the preparation of the order batches. The model is validated and tested with a real-life case study prepared from data acquired from the outpatient pharmacy of the Sint Maartenskliniek (SMK) in Nijmegen, The Netherlands, called the Maartensapotheek (MA). We also provide scenario analyses on the results, as well as two model variations to explore the influence of tactical-level decisions on the pharmacy’s profitability. Considering the above-mentioned significant characteristics of PSCN, these main research questions are as follows. 
\begin{itemize}
    \item How can medication synchronization be incorporated into the PSCN model?
    \item What delivery mode would be efficient and appropriate for various patient types?
    \item How can the medication needs of each patient type be considered in the model?
    \item How can each patient be effectively clustered based on their medication needs?
    \item What is the optimal number and composition for each batch of medication?
    \item How can staffing requirements for order assembling, preparation, and distribution be addressed? 
    \item What insights are crucial for the decision-makers and managers at Maartensapotheek (MA) to enhance their PSCN? 
\end{itemize}

\begin{table}[h]
\fontsize{6pt}{6pt}\selectfont
\ra{1.5}
\caption{\textcolor{black}{\label{PSCNLR} A Review on the recent PSCN studies.}}
\begin{tabular}{lcccccccccccccccc}
\hline
\multicolumn{1}{c}{\multirow{2}{*}{Study}} & \multirow{2}{*}{Formulation} & \multirow{2}{*}{} & \multicolumn{3}{c}{Decisions} & \multirow{2}{*}{} & \multirow{2}{*}{\begin{tabular}[c]{@{}c@{}}Medic.\\ Synch.\end{tabular}} & \multirow{2}{*}{} & \multicolumn{4}{c}{Configuration} & & \multicolumn{3}{c}{Solution Approach} \\ \cline{4-6} \cline{10-13} \cline{15-17} 
\multicolumn{1}{c}{} & & & AL & TR & LM & & & & Period & \begin{tabular}[c]{@{}c@{}}Deliv.\\ Mode\end{tabular} & \begin{tabular}[c]{@{}c@{}}Cooling\\ Mode\end{tabular} & \begin{tabular}[c]{@{}c@{}}Med. or\\ Product\end{tabular} & & ES & HR & MH \\ \hline
\cite{sousa2005global} & MILP & & * & * & - & & - & & M & - & - & - & & - & - & * \\
\cite{susarla2012integrated} & MINLP & & - & - & - & & - & & S & - & - & - & & - & * & - \\
\cite{hansen2015planning} & MILP & & - & - & - & & - & & S & - & - & - & & * & - & - \\
\cite{zhu2018design} & MILP & & - & * & - & & - & & M & - & - & - & & * & - & - \\
\cite{roshan2019two} & MILP & & * & * & - & & - & & M & - & - & - & & * & - & - \\
\cite{halim2019decision} & MILP/QD & & - & * & - & & - & & - & - & - & - & & * & - & - \\
\cite{goodarzian2020multi} & MINLP & & * & * & - & & - & & M & M & - & M & & - & - & * \\
\cite{nasrollahi2021mathematical} & MILP & & * & * & - & & - & & M & - & - & M & & - & - & * \\
\cite{ahmad2022multi}& MILP & & * & * & - & & - & & S & - & - & S & & * & - & - \\
\cite{li2023two} & MILP & & * & * & - & & - & & S & - & - & M & & * & - & - \\
\cite{abdolazimi2023development} & MILP & & * & * & - & & - & & M & - & - & M & & * & - & - \\
This Study & MILP & & * & * & * & & * & & M & M & M & M & & * & - & - \\ \hline
\multicolumn{17}{l}{\textit{\begin{tabular}[c]{@{}l@{}}Model Formulation: Mixed Integer Linear (MILP); Mixed Integer Non-Linear (MINLP); Qualitative/Decision-Making (QD); Allocation (AL), \\ 
Transportation (TR), Last-mile (LM); Single (S); Multi-(M)\end{tabular}}}                                        \\ \hline
\end{tabular}
\end{table}

To sum it up, the main contributions of this study based on the stated research questions can be listed as: (i) provide a mathematical formulation to manage the pharmaceutical cold supply chain while explicitly considering medication synchronization and different delivery modes, (ii) prepare, process, and provide a real case study of an outpatient pharmacy to be used for PCSCM-MSDM using data analytics methods, (iii) provide the exact optimal solutions for the understudied problem, (iv) perform extensive sensitivity (scenario) analysis, including variations on the base model, to derive insights for managers including but not limited to the outpatient pharmacy managers.
The remainder of this paper is structured as follows. Section 2 discusses relevant literature related to the model, application, and analytic methods presented in this paper. Section 3 presents the MILP model for the PCSCM-MSDM problem. Section 4 contains information about the MA case study, the data used, results and experiments, and alternative modeling choices. Finally, Section 5 concludes the paper, examines opportunities for further research, and discusses the limitations of this research.

\section{Problem Statement}\label{sec:PS}
This section explains our modeling approach for the Pharmaceutical Cold Supply Chain Management Considering Medication Synchronization and Different Delivery Modes (PCSCM-MSDM) problem. We propose a Mixed Integer Linear Programming (MILP) approach with the purpose of maximizing the logistical financial outcome (LFO) of an outpatient pharmacy considering various patient types, medication needs, the possibility of synchronizing medication, and different delivery modes. We also consider the staffing requirements for handling and distribution of orders in the optimization problem. 

Section \ref{subsec:description} explains the model and some background of the choices made. Section \ref{subsec:assumptions} states the assumptions made in the MILP model. Finally, Section \ref{subsec: mathematicalformulation} contains the mathematical formulation of the model: the sets, input parameters, decision variables, the objective function, and constraints.

\subsection{Model Description} \label{subsec:description}
The PCSCM-MDSM problem is inspired by the real application of an outpatient pharmacy, which serves a number of patients over a specified time horizon. The time horizon consists of multiple time intervals (we call them time periods) in which an order may be placed. The objective of the model is to optimize the logistical financial outcome of an (outpatient) pharmacy. The logistical financial outcome consists of the transportation cost, the handling cost, and the logistical turnover. The optimal solution yields the optimal number of order batches, the composition of these batches, and the number of employees needed to assemble and prepare the batches. The model contains patients who have a set of medications that they need to receive over the time horizon. Therefore, a patient orders all medications out of their set at least once over the time horizon. An order batch consists of all or a set of medications a patient needs repeatedly. From these sets of medications, an order batch can be composed, and the pharmacy receives a \textit{prescription line fee} for each medication type in a batch.
The PCSCM-MDSM problem is formulated based on the regulations and requirements of the Dutch healthcare system, specifically in relation to the logistical turnover of medications in pharmacies. This legislation allows pharmacies to ask for a prescription line fee from their patients to cover the logistical costs of providing medication \citep{Rijksoverheid2023WatRecept, ZorginstituutNederland2023ToelichtingReceptregelvergoeding}. The prescription line fee is issued per medication type for each ordered batch. The dosage of medication within an order is irrelevant to the prescription line fee. Hence, in our model, the prescription line fee only depends on the number of different medication types per order.

\begin{figure}[h!]
    \centering
    \resizebox{0.7\textwidth}{!}{
    \includegraphics[width = \textwidth]{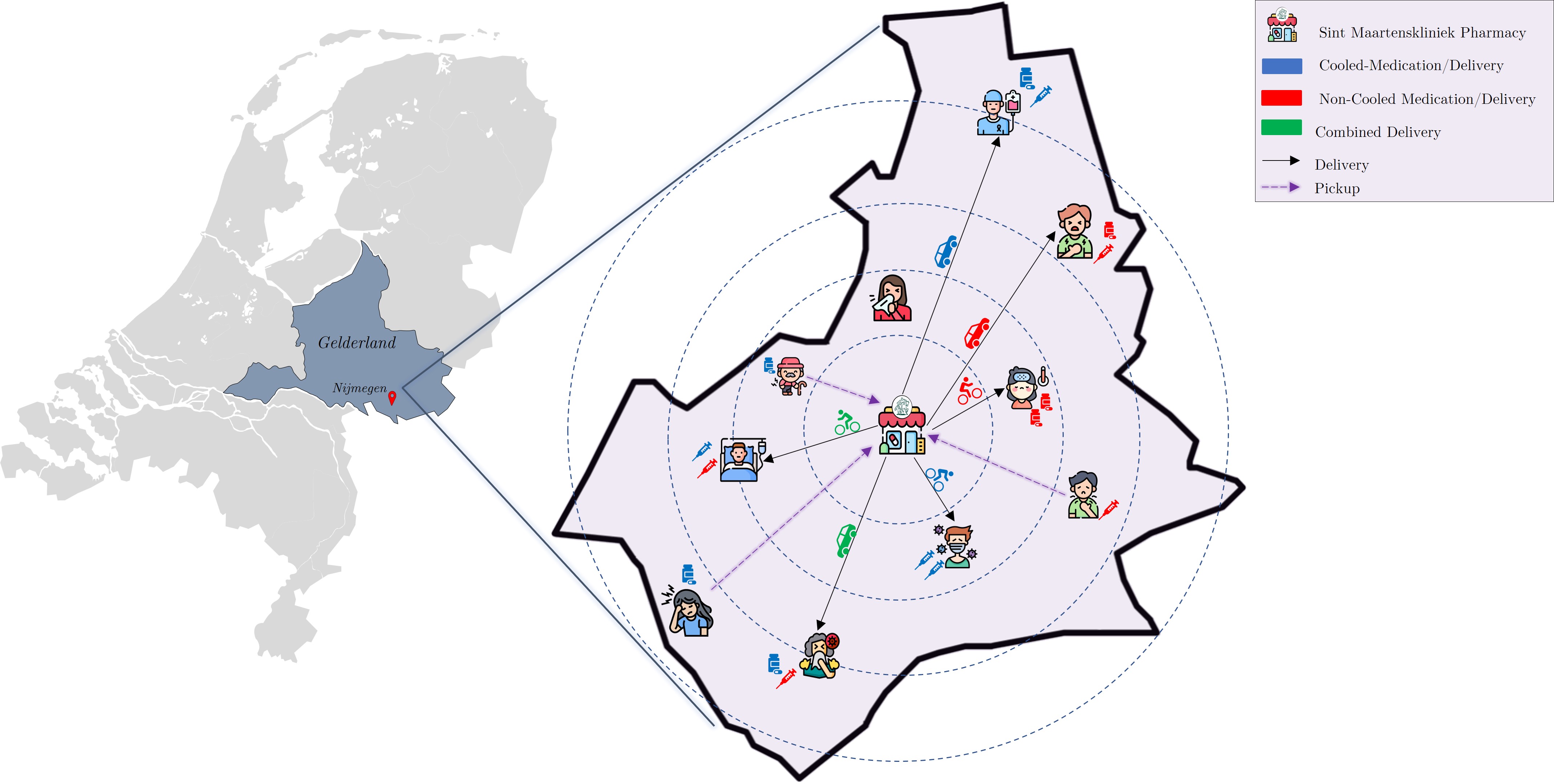}
    }
    \caption{\textcolor{black}{The schematic view of the proposed network.}}
    \label{fig:network}
\end{figure}

Moreover, medications have cooling requirements independent of the prescription line fee type. Therefore, we classify medications based on their cooling types (e.g., cooled or non-cooled). In addition, there is another classification required for medicine deliveries (e.g., cooled, non-cooled, combined). Medications need to be transported with a transportation cooling method fitting to their medication cooling type. The transportation cooling type is combined when a batch contains medications with different cooling requirements. \textcolor{black}{Such clustering of medications and transportation methods using temperature requirements aligns with the pharmaceutical guidelines currently in place for our case study. For future cases, we suggest developing a specific clustering approach, adapted to reflect the pharmacy's real-world conditions, concerning temperature-sensitive storage and transportation of medicines for different types of products and pharmacy guidelines.} The transportation costs depend on the delivery type (e.g., bike, car) and delivery cooling type of an ordered batch. \textcolor{black}{We note our pharmacy (like many pharmacies and hospitals) has a contract with an external carrier to ship the products. In such a deal, the shipping costs for different types of orders and delivery modes are agreed upon. Therefore, the transportation cost incurred by the pharmacy is not dependent on the distance covered; instead, it is hidden in the carrier's calculation of the fixed price for the transportation of medications accepted by the pharmacy.} Lastly, our model considers the number of hours needed for handling each ordered batch, which may differ depending on the transportation cooling type of the batch. As mentioned, we define the \textit{time period} as the shortest length of time that a patient may put one order. The time period should be in line with the practical requirements (e.g., a week, a month, etc.). The \textit{time horizon} is the total running time of the model and consists of at least one time period. \textcolor{black}{Finally, we can demonstrate the visual representation of the proposed PCSCM-MDSM problem as Figure \ref{fig:network}.}

We proceed by the assumptions considered to model PCSCM-MSDM, and the mathematical model in Section \ref{subsec: mathematicalformulation}, which includes notation, parameters, decision variables, and the explanation of the model.

\subsection{Assumptions} \label{subsec:assumptions}
In this section, we provide the assumptions considered in the model. We note that the assumptions are in line with the case study requirements that will be presented in Section \ref{sec:results}.
\begin{enumerate}
    \item The demand can always be satisfied from the pharmacy's stock.
    \item Our model is focused on satisfying demand while optimizing logistical revenues and costs, including handling, staffing, distribution, and delivery operations. Hence, the direct costs and earnings for medications (e.g., purchasing costs and health insurance reimbursements intended for the \textbf{direct} costs of medication) lie outside of the scope of this model.
    \item In practice, the size of a batch hardly affects the handling time of the batch, and the size also does not affect the delivery costs or the prescription line fee. Hence, we do not consider the size of a batch in the model (i.e., the number of doses, boxes, packages, and so on). Only the unique types of medications are considered in our model. We note that we consider the maximum number of batches that can be delivered through various delivery options in our model.
    \item Employees are assumed to be permanent staff. The number of employees of each type, therefore, needs to be the same in each time period.
    \item According to the case study data,  there has been no documented occurrence of harm or loss to medications during the delivery process. Therefore, the model does not incorporate any cost components related to such incidents.
    \item Like many other pharmacies and hospitals, an external carrier is hired to transport medicines. Hence, the model does not include routing decisions to be optimized. 
\end{enumerate}

\subsection{Mathematical Formulation} \label{subsec: mathematicalformulation}
This section presents the mathematical formulation of the PCSCM-MDSM problem by providing the sets, indices, decision variables, the objective function, and constraints.

First, we introduce the sets and indices used in the model that can be found in Table \ref{table:sets}.

\begin{table}[!h] 
\caption{Sets and indices of the MILP model}
\label{table:sets}
\scalebox{0.7}{
\begin{tabular}{ll}
\toprule
Set & Description                                          \\ \midrule \addlinespace[0.05cm]
 $A$   & set of  transportation cooling types, $a \in A$       \\
$C$   & set of cooling types related to medications, $c \in C$  \\
$D$   & set of delivery types, $d\in D$                        \\
$E$   & set of employee types, $e\in E$                        \\
$W$   & set of time periods, $w\in W$                          \\
$K$   & set of medication types, $k\in K$                      \\
$P$   & set of patient types, $p\in P$                         \\
\bottomrule
\end{tabular}}
\end{table}



Here, we explain the definition and notations of the input parameters.
The time horizon for which we solve the model consists of $W$ time periods. Parameter $q_{ckp}$ states the number of unique medications of cooling type $c$ and medication type $k$ that a patient type $p$ uses consistently over the time horizon. Hence, $\sum_{c\in C} \sum_{k\in K} q_{ckp}$ is the total number of unique medications that a patient $p$ needs over the time horizon.
We note that the pharmacy planning time horizon for special medications is relevant as it is in line with their patient ordering policy and a significant number of products share the same consumption period (For our case it is typically spanning between 4 or 6 months). 

Our model has three revenue/cost components with associated parameters. The prescription line fee is revenue of the pharmacy that is collected from the insurer, patient, or hospital. We denote this fee with $f_k$, which depends on the number of unique medication types per batch. Hence, when the number of unique medication types increases in a batch, this fee also increases. The quantity of each unique medication is irrelevant, as the prescription line fee is meant to cover administrative costs regarding the processing of a certain type of medication. 
 The transportation cost, $t_{ad}$, is incurred in each transportation of a batch from the pharmacy to the home of the patient. Therefore, $t_{ad}$ depends on the combination of transportation cooling type $a$ and transportation method $d$. The transportation cooling type $a$ relates to the type of cooling assigned to the transportation of the batch and the transportation method $d$ relates to the type of transportation, such as truck, bicycle, or drone delivery. 
 Finally, $s_{e}$ is the staffing cost of employee type $e$ time horizon. 



Here, we define the required upper bounds and lower bounds for the parameters in the model. First, our preliminary data analysis shows that some patients need similar types of medications per time period; therefore, we classified similar patients into one group. This approach helps to reduce the number of decision variables and, in turn, the problem instance size. Hence, we define parameter $\rho_{p}$ as the number of patients within patient type $p$. 
Depending on the patient type, we consider a lower bound on the number of orders per time period. The parameter values are extracted based on the behavior of different types of patients. For instance, some types of patients may forget part of their medications in one order, or some may receive specific test medications; therefore, $\sigma_{p} \ge 1$ allows more than one order over the time horizon. Second, there is a maximum level $\delta_{dw}$ for the number of batches that can be delivered with transportation method $d$, in time period $w$. Third, the maximum number of hours an employee $e$ can work per period $w$ is $\theta_{ew}$. We note that it takes $u_{ae}$ hours for employee type $e$ to process one batch of type $a$ considering the communication with the patient, back office work, packaging, and auxiliary activities. We summarize the description of input parameters in Table \ref{table: parameters}.

\begin{table}[!h]
\ra{0.85}
    \caption{
    Input parameters of the MILP model}
     \label{table: parameters}
     \scalebox{0.8}{
    \begin{tabular}{lll}
    \toprule
    Parameter                        & Description                                                                        \\ \midrule \addlinespace[0.05cm]
    $q_{ckp}$     & Number of unique medicines with cooling type $c$, medication                              
                                 type k, and for patient type $p$ \\
    $f_{k}$                      & Prescription line fee for medication type $k$                                      \\
    $\rho_{p}$                      & Number of patients of type $p$                                                     \\
    $\sigma_{p}$       & Minimum number of orders per time horizon of patient type $p$                        \\
    $t_{ad}$           & Transportation cost with transportation cooling $a$ and delivery type $d$   \\
    $\delta_{dw}$     & Maximum number of orders to be delivered with method $d$ 
                                                   in time period $w$          \\
    $u_{ae}$         & Time (in hours) needed of employee type $e$ to process one order 
                                                   of type $a$       \\
    $\theta_{ew}$    & Maximum number of hours an employee of type $e$ can work 
                                                   per time period $w$            \\
    $s_{e}$          & Salary of employee type $e$ for the time horizon                                             \\
    \bottomrule
    \end{tabular}}
\end{table}

The proposed model has four sets of decision variables. The first set of decisions pertains to whether an order is placed for patient $p$ in period $w$, considering the transportation cooling type $a$ and transportation method $d$. This decision is represented by the variable $x_{adpw}$. If the order is placed, the value of $x_{adpw}$ is set to one; otherwise, it equals zero. The decision regarding the composition of the orders is denoted by variables $o_{ckpw}$. We define $o_{ckpw}$ as the number of unique medications of cooling type $c$ and type $k$ that are ordered by patient $p$ in period $w$. The third decision encompasses the number of employees needed to handle and distribute the batches per time period. This decision is modeled by the variables $m_{ew}$, which is the number of employees of type $e$ needed in time period $w$. The final decision variable set is the number of employees needed in a cycle. This decision is modeled by variables $M_e$, which is the number of employees $e$ needed in a cycle. The decision variables for this model can be found in Table \ref{table: decvar}.



\begin{table}[H]
\ra{1.1}
    \caption{
    Decision variables of the MILP model}
    \label{table: decvar}\
    \scalebox{0.65}{
    \begin{tabular}{lll}
        \toprule
        Binary Variable                  & Description \\ \midrule \addlinespace[0.05cm]
        $x_{adpw}$    &   1 if an order is placed in period $w$ for patient type $p$ with 
                          delivery type $d$, and transportation cooling $a$, 0 otherwise \\ 
                    \midrule \addlinespace[0.05cm]
        Integer Variables                     & Description \\  \midrule \addlinespace[0.05cm]
        $o_{ckpw}$  &  Number of unique medicines in order, with cooling type   
                      $c$, type $k$, ordered at time $w$, and for patient type $p$. \\
        $m_{ew}$    &  Number of employees of type $e$ needed in period $w$.\\
        $M_{e}$     &  Number of employees of type $e$ needed in the time horizon. \\ \bottomrule
    \end{tabular}}
\end{table}

The objective of the model is to maximize the LFO of an outpatient pharmacy that includes one revenue and two cost parts. The revenue part is the prescription line fee per unique medication per order (the first part of the objective function). The cost part includes the transportation costs per order (second part of the objective function) and the employee costs related to the distribution of medication (third part of the objective function. Equation \eqref{ObjFunc} shows the objective function.

\begin{align}
\text{max} &\sum_{c\in C} \sum_{k\in K}\sum_{p\in P}\sum_{w\in W} o_{ckpw} f_{k} \rho_{p}- \sum_{a\in A} \sum_{d\in D} \sum_{p\in P}\sum_{w\in W}  x_{adpw}t_{ad}\rho_{p}- \sum_{e \in E}M_{e}s_{e}\label{ObjFunc}
\end{align}

We have defined nine constraint sets to ensure that the model behaves as desired. Constraint set \eqref{COneDelivery} limits the number of batches to one per time period for each patient, and only one combination of transportation cooling $a$ and delivery type $d$ can be selected for each batch. Constraint sets \eqref{CAandC1} and \eqref{CAandC2} ensure that the number of unique medications in a batch, $o_{ckpw}$, relates to $x_{adpw}$. That means, if no order is placed in period $w$ for patient $p$, $x_{adpw}$ equals 0, then $o_{ckpw}$ must be zero since there cannot be any medications in an order that does not exist. Also, Constraints \eqref{CAandC1} and \eqref{CAandC2} make sure that the right type of transportation cooling $a$ is related to the right type of medication cooling $c$ in a batch. To clarify, a batch with only cooled medicine should be transported with cooled packaging, an order with only non-cooled medicine should be transported with non-cooled packaging, and an order with both cooled and non-cooled medication should be transported with combined packaging. Furthermore, constraint set \eqref{CUniqueMeds1} ensures that the medication requests for each patient type $p$ for every medication type $k$ with cooling type $c$, $q_{ckp}$, are met over the entire time horizon, ensuring that each patient receives the medications they need.
We define the level of medication synchronization, $\sigma_{p}$, as the number of times that patient type $p$ needs to receive batches of orders; therefore, Constraint \eqref{CMinNumOrders} ensures that the total number of order batches for each patient type $p$ equals or is greater than $\sigma_{p}$. 
The constraint ensures that the patients receive batches multiple times over the time horizon. Due to the limited capacity of the pharmacy to deliver medicines using different delivery modes, we use Constraint \eqref{CDeliveryMethods} to ensure that there are no more than $\delta_{dw}$ batches delivered with type $d$ in each period $w$.
Constraint \eqref{CHoursPerWorker} makes sure that for every period $w$, there are enough employees available to perform the labor related to packaging and distributing the orders.
Constraint \eqref{CTotalWorkers} makes sure that the number of employees of type $e$ in the time horizon is greater than or equal to the number of employees of type $e$ needed in each period $w$. 
Constraints \eqref{SignConstraints} contain the type and sign restrictions of this model.



\begin{align}
    \text{s.t.} \quad & \sum_{a\in A}\sum_{d\in D} x_{adpw} \le 1  & \quad &  \forall p \in P, w \in W \label{COneDelivery}\\
    & \sum_{c\in (c1)}\sum_{k\in K} o_{ckpw} \le  M_1 \sum_{a\in (a1, a3)}\sum_{d\in D} x_{adpw}      & \quad & \forall p\in P, w \in W \label{CAandC1}\\
    & \sum_{c\in (c2)}\sum_{k\in K} o_{ckpw} \le M_2 \sum_{a\in (a2, a3)}\sum_{d\in D} x_{adpw}  & \quad &\forall p\in P, w \in W \label{CAandC2}\\
    & \sum_{w\in W}o_{ckpw} = q_{ckp},   & \quad & \forall c\in C, k\in K, p\in P \label{CUniqueMeds1}\\
    & \sum_{a\in A}\sum_{d\in D}\sum_{w \in W} x_{adpw} \ge \sigma_{p},                       &\quad &  \forall p\in P \label{CMinNumOrders}\\
    & \sum_{a\in A}\sum_{p\in P} x_{adpw} \rho_{p} \le \delta_{dw}                    & \quad & \forall d\in D, w\in W \label{CDeliveryMethods}\\
    & \sum_{a\in A}\sum_{d\in D}\sum_{p\in P} x_{adpw}u_{ae}\rho_{p} \le m_{ew}\theta_{ew}    &\quad &  \forall w \in W \label{CHoursPerWorker}\\
    & M_{e} \ge m_{ew}           & \quad & \forall e \in E, w\in W \label{CTotalWorkers}\\
    & x_{adpw} \in \{0,1\}, o_{ckpw} \in \mathbb{Z}^{+}, m_{ew} \in \mathbb{Z}^{+}, M_{e} \in  \mathbb{Z}^{+} \label{SignConstraints}
\end{align}

\section{Case Study and Results} \label{sec:results}
In this section, we present the case study of the Maartensapotheek (MA) outpatient pharmacy and optimize their LFO using the MILP model explained in Section \ref{sec:PS}. The MA pharmacy of the Sint Maartenskliniek (SMK) is located in Nijmegen, the Netherlands. The SMK is a hospital specializing in mobility and posture and has three main departments: orthopedics, rheumatology, and rehabilitation. The MA provides medication for patients of the SMK. This mostly concerns rheumatology patients, as the rheumatology department is the largest outpatient department of the SMK. Besides, medication for patients of the orthopedic and revalidation departments is often provided via the inpatient pharmacy of the SMK. This section explains the problem requirements of the MA, the data used for the case study, and the model scenarios. 

The ultimate goal of MA is to provide a high-quality service to its patients by giving them the option of at-home delivery of medications for free. However, as this service is deemed to be expensive; therefore, the pharmacy wishes to get insight into its logistical costs and revenues to optimize the LFO. Hence, we propose the PCSCM-MDSM problem formulation to solve the MA case. 

Section \ref{subsec:data} outlines the input parameters used for the MA case and explains the data preparation procedure. In Section \ref{subsec: scenarios}, we describe the scenarios used in the case study and explain the need for scenario analysis. Section \ref{subsec:results} presents the results of all scenarios and relevant metrics. To explore the potential for future development of this case or its application to similar cases, Section \ref{subsec:modelvariations} introduces two variations that slightly modify the model settings and relax certain assumptions. We compare the results of these model variations to the original results and provide managerial insights. Finally, in Section \ref{subsec:methods statistics}, we report on the applied solution methods and provide insights using reported statistics. We proceed with explaining the data used in our model.

\subsection{Data} \label{subsec:data}
This section explains how and what data is collected and processed for the MA case. First, in the data collection phase, we explain our method and sources to gather relevant data. The data is mostly used for determining the values of input parameters that come from the medication orders of MA in 2021 and 2022. The data list of orders contains order batches from patients, which include at least one medication product per batch. The data contains three types of unique identifiers as follows.

\begin{itemize}
    \item  Order number: It is used to derive the order date, the delivery type, and the number of orders.
    \item  Patient number: It is used to derive the location of the patient (which can be relevant for certain delivery types), the number of unique patients, the medications that the patients need, and the order composition per patient (in combination with the medication product numbers in the orders of each patient). 
    \item Medication product: It is used to find the medication type and the medication cooling restriction. 
\end{itemize}

The set elements used in the model are derived from patient and order data, observations of the order preparation and delivery process, and discussions with expert employees of the MA. Sets $A$ and $C$ are derived from the cooling restrictions and transportation cooling types that are used in the MA. Set $D$ consists of the delivery types the management of the MA considers. Set $E$ contains the two types of employees present in the MA. Set $K$ is based on the existence of the prescription line fee of a medication. The MA has two types: medication with and without a prescription line fee. Set $W$ consists of four months, as this is the current maximum order period in the MA. The sets and their elements for the MA case are shown in Table \ref{table: MACaseSets}. We explain set P further in the explanation of the parameter $q_{ckp}$, $\rho_p$, and $\sigma_p$.

\begin{table}[h!]
\ra{0.85}
    \caption{Set elements for the MA case}
    \label{table: MACaseSets}
    \scalebox{0.8}{
    \begin{tabular}{ll}
        \toprule
        Set   & Elements   \\
        \midrule \addlinespace[0.05cm]
        $A$ & \{Cooled  transportation, Non-cooled transportation, Combination transportation\}                      \\
        $C$ & \{Cooled medication, Non-cooled medication\}                                                           \\
        $D$ & \{Truck transportation, delivery hubs, bicycle transportation, pick-up\}                             \\
        $E$ & \{Pharmaceutical employee, pharmacy technician\}                                                      \\
        $K$ & \{Medication w/o prescription line fee, Medication with  prescription line fee\}                  \\
        $P$ & \{0, 1, \dots, 224\} \\                                                           
        $W$ & \{January, February, March, April\} \\  
        \bottomrule
    \end{tabular}}
\end{table}

The consumption period of the medications is between 4 to 6 months. Other pharmacies can anticipate their own planning horizons that could be as short as one month or even one week, allowing for consideration of varying consumption periods for different medications. Every patient $p$ has a value of $q_{ckp}$ for each $c$ and each $k$, which is shown in Table \ref{tab: summaryqrhosigma}. We derive these values by first categorizing patients into four groups based on the medications they need. As there are two elements in sets $C$ and $K$, the total number of combinations is four. Second, we consider how frequently patients have placed orders in the past, which relates to $\sigma_p$. In this regard, we use aggregation to create patient types such that irregular or infrequent orders are grouped with the closest patient types to ensure that set $P$ is not too large. For example, a patient type could be a patient who needs two medications of $c_0$ and $k_0$ and three medications of $c_1$ and $k_0$ that typically orders twice every four months. From the data, we see that there are five patients with the same specifications and one patient who has the same medication types but orders three times every four months. We place that last patient into the same category $p \in P$, as we see this patient as an outlier. This means that $\rho_p$ is 5 + 1 = 6. Using this approach, we identify 225 patient types $p$. We determine $\sigma_p$ by determining the number of times each patient $p$ typically orders over four months period, and $\rho_p$ by the number of patients in the same index $p$. Table \ref{tab: summaryqrhosigma} presents only a portion of values related to $q_{ckp}$, $\sigma_{p}$, and $\rho_{p}$ to enhance readability. Hence, as we only show parts of the rows, it is coincidental to have one as the minimum number of orders (shown in the rightmost column) in each row. Although not explicitly displayed in the table, it should be noted that $\sigma_p = 1$ for 157 patient types, $\sigma_p = 2$ for 57 patient types, and $\sigma_p = 3$ for 11 patient types. The comprehensive sets of values can be found in Tables A1 and A2 in the Appendix.



\begin{table}[h!]
\caption{Summary of $q_{ckp}$, $\rho_p$ and $\sigma_p$: the number of unique medicines with cooling type $c$, medication type $k$ and for patient type $p$; the number of patients of type $p$ and the minimum number of orders for patient $p$.}
\label{tab: summaryqrhosigma}
\scalebox{0.8}{
\begin{tabular}{lllllcc}
\toprule
  \textbf{Patient type $\boldsymbol{p}$}                                             & \multicolumn{4}{l}{\textbf{Nr. unique meds $\boldsymbol{q_{ckp}}$}}                                       & \textbf{Nr. patients $\boldsymbol{\rho_p}$} & \textbf{Min. nr. orders $\boldsymbol{\sigma_p}$} \\ \cline{2-5} 
     & $c_0 ; k_0$       & $c_0 ; k_1$  & $c_1 ; k_0$  & $c_1 ; k_1$   &&            \\ \midrule
$p_{0}$                                      & 1                              & 0                         & 0                         & 0                         &935&1\\
$p_{1}$                                      & 1                              & 0                         & 0                         & 0                         &553&1\\
$p_{2}$                                      & 0                              & 0                         & 0                         & 1                         &539&1\\
$p_{3}$                                      & 0                              & 0                         & 0                         & 1                         &318&1\\
$p_{4}$                                      & 0                              & 0                         & 0                         & 2                         &261&1\\
\dots                                        & \dots                          & \dots                     & \dots                     & \dots                      & \dots & \dots             \\
$p_{224}$                                    & 0                              & 1                         & 0                         & 6                          &1&1 \\
\bottomrule
\end{tabular}}
\end{table}



The prescription line fee is determined by grouping the medications into two sets: medications without a prescription line fee ($k_0$) and medications with a prescription line fee ($k_1$). Obviously, $f_{k_0}$ equals zero. Within the medications with a prescription line fee $k_1$, there are differences that are often based on the type of order or other factors that are difficult to trace. Based on the available data, we estimate the differences in prescription line fees are quite low (coefficient of variation of 0.45). Therefore, we use a weighted average of historical values of the prescription line fee and find $f_{k_1} = $ \geneuro 7.94. 

We obtain the transportation costs, $t_{ad}$, by using the agreements made with the delivery services employed by the MA. The pick-up costs are estimated by multiplying the hourly salary of a pharmaceutical employee with the time it takes for a pick-up, which is on average 3 minutes. It is important to note that a pick-up at the pharmacy is not the same as an at-home delivery or hub pick-up, which takes only a few seconds. This is because the pharmacy chooses to provide the additional service of verifying the order and patient details with individuals at the time of medication pick-up from the MA. The maximum number of transportation (travel between pharmacy and drop-off points) with each transportation method $d$, $\delta_{dw}$, is determined based on agreements with the delivery services, which is typically fixed across all periods, $w$. In determining the monthly maximum number of transportation, we consider the locations of the patients and assess whether bicycle delivery and pick-up from the hubs are feasible options for those patients. Pick-up is not feasible for distances greater than 10 kilometers from the pick-up point (i.e., the MA location) or pick-up hub. Furthermore, bike delivery is not possible for distances exceeding 15 kilometers from the bicycle starting point (i.e., the MA location). The limit for truck transportation is much higher, and it is practically infinite, but we set it to a sufficiently large value to ensure that all transportation can always be conducted with this method. The values for $t_{ad}$ and $\delta_{dw}$ are provided in Tables \eqref{tab: transportationcost_t} and \eqref{tab: delta}.


\begin{table}[h!]
\centering
\begin{minipage}{0.46\linewidth}
\ra{0.85}
    \caption{Values of $t_{ad}$: the transportation cost with transport cooling $a$ and delivery type $d$}
    \label{tab: transportationcost_t}
    \scalebox{0.55}{
    \begin{tabular}{llll}
        \toprule
        & \multicolumn{3}{l}{\textbf{Transportation cooling type $\boldsymbol{a}$}}\\ \cline{2-4} \addlinespace[0.05cm]  
        \textbf{Delivery type $\boldsymbol{d}$}       & Cooled ($a_0$)& Non-cooled ($a_1$) & Combination ($a_2$)           \\ \midrule 
        Truck transportation ($d_{0}$)            & \euro \space  16.62   & \euro \space  11.64   & \euro \space17.32                     \\
        Delivery hubs ($d_{1}$)              & \euro  \space 3.00   & \euro \space 3.00    & \euro \space3.01                      \\
        Bicycle transportation ($d_{2}$)            & \euro \space 13.62   & \euro \space  8.64    & \euro \space 14.32                     \\
        Pick-up ($d_{3}$)              & \euro \space 1.65    & \euro \space 1.65    & \euro \space 1.66                     \\
        \bottomrule
    \end{tabular}
    }
\end{minipage}
\hfill
\begin{minipage}{0.49\linewidth}
\ra{1.05}
    \caption{Values of $\delta_{dw}$: the max number of orders to be delivered with method $d$ in time period $w$}
    \label{tab: delta}
    \scalebox{0.65}{
    \begin{tabular}{lc}
        \toprule
        \textbf{Delivery type $\boldsymbol{d}$}       & \textbf{$\boldsymbol{\delta_{dw}}$ ($\boldsymbol{\forall w \in W}$)}             \\ \addlinespace[0.05cm]  \midrule
        Truck transportation ($d_{0}$)     & 6950                     \\
        Delivery hubs ($d_{1}$)            & 271                      \\
        Bicycle transportation ($d_{2}$)   & 120                     \\
        Pick-up ($d_{3}$)                  & 253                     \\
        \bottomrule
    \end{tabular}
    }
\end{minipage}
\end{table}


We determine the time required for employee type $e$ to process a batch with cooling type $a$, denoted as $u_{ae}$, by measuring the order preparation time of all activities. Moreover, we determine $\theta_{ew}$ from the actual working hours of each employee type per period, $w$. Finally, we determine $s_e$ based on each employee type's total hourly salary costs. The values for $u_{ae}$, $\theta_{ew}$, and $s_{e}$ are presented in Table \ref{tab:Real}.

\begin{table}[h!]\centering
\ra{0.85}
\caption{The time (in hours) needed of employee $e$ to process one order of type $a$; the maximum number of hours an employee of type $e$ can work per time period $w$; the hourly salary per employee type $e$ }
\scalebox{0.65}{
\begin{tabular}{@{}c>{\raggedleft}p{2cm}>{\raggedleft}p{2cm}>{\raggedleft}p{2cm}c>{\raggedleft}p{4cm}>{\raggedleft}p{1cm}c>{\raggedleft}p{2cm} @{}}\toprule
 & \multicolumn{3}{c}{\textbf{Hours needed to process order $\boldsymbol{a}$ by employee $\boldsymbol{e}$, $\boldsymbol{u_{ae}}$}} & & & &   \\ \cline{2-4}	 \addlinespace[0.1cm]  
	 \textbf{Employee type $\boldsymbol{e}$} 	&	Cooled ($a_0$)& Non-cooled ($a_1$) & Combination ($a_2$)	& &	\textbf{Maximum number of hours $\boldsymbol{\theta_{ew}}$ ($\boldsymbol{\forall w \in W}$)} &  & \textbf{Salary $\boldsymbol{s_{e}}$}\\	\midrule \addlinespace[0.1cm]   \addlinespace[0.05cm]
        Pharmaceutical employee ($e_{0}$) &  0.1293   & 0.1477    & 0.1779 & & 126.667 & &    \euro \space 33.00  \\
        Pharmacy technician ($e_{1}$)     &  0.0233   & 0.0233    & 0.0233 & & 126.667 & &    \euro \space 40.0 \\	
\bottomrule
\end{tabular}}
\label{tab:Real}
\end{table}

Finally, the values for $M_1$ and $M_2$ are both set to 15. This is because $\sum_{c \in C}\sum_{k\in K} q_{ckp}$ out of the data has a maximum value of 15 $\forall p \in P$. Therefore, the values for $M_1$ and $M_2$ can be set to 15.

\subsection{Computational Results} \label{subsec:results}
In this section, we present the results of the model using the input data described in Section \ref{subsec:data}, referred to as the `base case'. Additional scenarios are discussed in Section \ref{subsec: scenarios}. We begin by introducing the key performance indicators (KPIs) used to evaluate the results. These KPIs are divided into two categories: output KPIs and decision KPIs. The output KPIs provide insights into the outcomes of implementing specific decisions, while the decision KPIs highlight the decisions required to achieve desired outcomes. The first output KPI is \textit{total annual LFO}, which is the LFO from the model multiplied by 3, as the model for the MA runs for four months. The total annual LFO consists of the annual transportation costs, the annual handling costs, and the annual prescription line fee. This KPI shows the pharmacy the total picture of its annual LFO. The second output KPI is \textit{LFO per order}, which is the LFO divided by the number of orders. This KPI shows the LFO relative to the orders, which can give a different perspective than the total LFO. The decision KPIs are \textit{annual number of order batches}, \textit{annual number of order batches per cooling type}, \textit{annual number of transportation per transportation type}, and \textit{number of employees per time period}. These KPIs provide insights into enhancing synchronization strategies, the composition of order batches, delivery modes, and staffing requirements. 

Table \ref{tab: base case valuesKPI12} presents the results for KPIs 1 and 2, which include the breakdown of LFO into transportation cost, handling cost, and prescription line fee (revenue). The results indicate a negative total LFO and negative LFO per order, indicating that the pharmacy is operating at a loss in this aspect. However, it's important to note that the pharmacy is part of a larger hospital where the logistics expenses are compensated for by other areas, mitigating the impact of the logistical loss on the pharmacy or hospital as a whole. Additionally, this paper solely focuses on logistical outcomes and does not consider results generated from procurement and sales. Furthermore, the total LFO in 2022 was \euro \space -197,208.23, indicating that the utilization of the model increases the MA's LFO by 34\%.


\begin{table}[h!]
\ra{1.05}
    \caption{Base case LFO for output KPIs: KPI 1 and 2}
    \label{tab: base case valuesKPI12}
    \scalebox{0.75}{
        \begin{tabular}{lll}
        \toprule
        \textbf{}             & \textbf{KPI 1: Total annual LFO} & \textbf{KPI 2: LFO per order} \\ \midrule \addlinespace[0.1cm]
        Total LFO          & \euro \space         -130,874.47 & \euro  \space          -5.76         \\
        Transportation cost   & \euro \space         -231,387.72 & \euro  \space          -10.18        \\
        Handling cost         & \euro \space         -204,835.33 & \euro  \space          -9.01         \\
        Prescription line fee & \euro \space \space  305,348.58  & \euro  \space  \space  13.43         \\ \bottomrule
        \end{tabular}}
\end{table}

Tables \eqref{tab: base case valuesKPI345} and \eqref{tab: base case valuesKPI6} show the results for the decision KPIs. To achieve the LFO in Table \ref{tab: base case valuesKPI12}, the pharmacy needs to have 22,740 orders on a yearly basis. The number of orders in 2022 was 21,347, which is close to the model output. This shows that the order frequency resulting from the model is attainable for the MA. Furthermore, we see that most of the batches are delivered with non-cooled transportation ($a_1$), which was also the case in 2022. The number of batches delivered with this cooling method decreased from 49.5\% to 46.4\%, while the number of batches transported with combination transportation ($a_2$) increased slightly, from 19.7\% to 23.8\%. This shows that the model combines more orders thus synchronizing order batches more efficiently. The number of cooled transportations ($a_0$) was hardly affected. As for the transportation methods, we see that the most used transportation method is truck transportation ($d_0$). This is the most expensive transportation mode, but it is the only unrestricted mode of transport. We see that the other transportation modes are maximized to their full capability, showing the potential of alternative, cheaper transportation modes. Finally, we note that there are 3 pharmaceutical employees needed and 1 pharmacy technician, for the logistical operations of the MA. This is the current occupation at the MA so the staffing does not need to change for this solution.

\begin{table}[h!]
\centering
\begin{minipage}{0.45\linewidth}
\caption{Base case results for decision KPIs 3, 4, and 5}
\label{tab: base case valuesKPI345}
\scalebox{0.6}{
\begin{tabular}{lll}
\toprule
& \textbf{Number of order batches} & \textbf{\% of orders} \\
\hline \addlinespace[0.05cm]
\multicolumn{3}{l}{\textbf{KPI 3: Annual number of order batches}}\\ 
Total nr. of orders                  &       22,740  &            100.0\% \\
\midrule \addlinespace[0.05cm]
\multicolumn{3}{l}{\textbf{KPI 4: Annual number of order batches per cooling type}}\\  \addlinespace[0.05cm]
Cooled transport ($a_0$)            &       6,783   &             29.8\% \\
Non-cooled transport ($a_1$)        &      10,551   &             46.4\% \\
Combination transport ($a_2$)       &       5,406   &             23.8\% \\
\midrule \addlinespace[0.05cm]
\multicolumn{3}{l}{\textbf{KPI 5: Annual number of transportations per transportation types}}\\  \addlinespace[0.05cm]
Truck transportation ($d_{0}$)      &       15,012  &            66.0\%         \\
Delivery hubs ($d_{1}$)             &       3,252   &            14.3\%        \\
Bicycle transportation ($d_{2}$)    &       1,440   &            6.3\%         \\
Pick-up ($d_{3}$)                   &       3,036    &           13.4\%       \\ \bottomrule
\end{tabular}}
\end{minipage}
\hfill
\begin{minipage}{0.44\linewidth}
\caption{Base case results for decision KPI 6}
\label{tab: base case valuesKPI6}
\scalebox{0.6}{
\begin{tabular}{lc}
\toprule
\textbf{Employee $\boldsymbol{e}$}             & \textbf{\# employees per time period} \\
\midrule \addlinespace[0.05cm]
Pharmaceutical employee ($e_{0}$)     &       3    \\
Pharmacy technician ($e_{1}$)          &       1    \\ \bottomrule
\end{tabular}}
\end{minipage}
\end{table}

The running time of this solution is 5.29 seconds, with an optimality gap of $1.83*10^{-15}$, showing that the solution is (virtually) optimal. We further analyze solution statistics, such as optimality gap, running time, and instance size, in Section \ref{subsec:methods statistics}.

Based on our analysis, we conclude that the MA can achieve an LFO increase of 34\% by implementing small adjustments, particularly by exploring alternative transportation methods. These changes primarily involve optimizing the number of order batches. In Section \ref{subsec: scenarios}, we evaluate more substantial modifications that the MA or other pharmacies can consider to enhance profitability and efficiency further.


\subsection{Scenario Analysis} \label{subsec: scenarios}
This section presents the scenario analysis, which involves 23 modifications to the data discussed in Section \ref{subsec:data}. The purpose of the scenario analysis is to assess input changes that require more significant investments compared to the minor adjustments made to the pharmacy logistics in the model results presented in Section \ref{subsec:results}. Additionally, the scenario analysis allows us to evaluate the robustness of the model. We categorize the adaptations into four groups: the inclusion of annual unique patients, synchronization level, patient types, and order period. We begin by explaining these four types of adaptations and their impact on specific input parameters. Subsequently, we analyze and discuss the effects of the scenarios on the model results.

The first modification involves changing the number of patients in the model. Currently, the annual number of unique patients is 7,000. However, there are plans for the hospital to take over a department from a nearby hospital, which is expected to increase the number of patients. The projected annual number of patients is 10,000. This change in the annual number of patients impacts the parameter $\rho_p$ in the model. Specifically, the current value of $\rho_p$ is extrapolated to accommodate the new patient count. The sum of all $\rho_p$ values, denoted as $\sum_{p \in P} \rho_p$, is adjusted from 7,000 to 10,000 to reflect the increased patient population. 

The second modification involves changing the synchronization level, which is divided into three cases: 77\%, 87\%, and 100\% synchronization. The synchronization level is determined by the completeness of an order, which is calculated based on the number of medications ordered compared to the total number of medications on the prescription. For example, if a patient has four medications prescribed but only orders three, the completeness of that order is 75\%, which corresponds to the synchronization level. The overall synchronization percentage of 77\% is a weighted average of these completeness percentages. The synchronization parameter, denoted as $\sigma_p$, is derived from this percentage, and patients are grouped into patient types, denoted as $p$, based on the synchronization level and the parameter $q_{ckp}$. To evaluate the potential for synchronization improvement, we adjust the $\sigma_p$ parameter for each patient type $p$ according to two scenarios: realistic improvement and ideal improvement. In the realistic improvement scenario, we collaborate with pharmacy management to assess feasible synchronization enhancements. For instance, a patient who currently orders three times in the time horizon may reduce their orders to two times. The average synchronization level for this realistic scenario is 87\%. In the ideal improvement scenario, all patients order only once throughout the time horizon, resulting in 100\% synchronization. To assess these adaptations, we adjust the parameter $\sigma_p$ based on the specified synchronization level.

Thirdly, we consider an adaptation in the patient types served by the pharmacy. Currently, the MA serves all patients from the SMK or related hospitals who wish to order medication at the MA. However, the hospital is exploring a scenario where they only cater to patients with medication type $k_0$, which refers to medication without a prescription line fee. These are medications that cannot be obtained from other pharmacies, meaning patients are required to order them exclusively from the MA. Patients with other medication types ($k_1$), which include medications with a prescription line fee, have the freedom to visit their local pharmacy if they prefer. This scenario has an impact on parameter $\rho_p$ as patients who do not have any medications of type $k_0$ in their prescription are excluded, resulting in their $\rho_p$ value being set to 0. Consequently, this adaptation decreases the total number of patients, leading to a reduction in the sum of $\rho_p$ values, denoted as $\sum_{p \in P} \rho_p$.

The final adaptation involves adjusting the length of the time horizon. Presently, the MA allows patients to order their medication for a maximum period of 4 months. However, there is an expectation that costs will decrease if the time horizon is extended to 6 months. This adaptation does not have any impact on the input parameters but rather expands the number of elements in set $W$ from 4 to 6.

A total of 24 scenarios are formed by combining these four adaptations: two options for the number of patients, three options for the synchronization level, two options for the patient types, and two options for the order period. Table \ref{table: scenarios} provides a summary of the scenarios, including the total annual LFO for each scenario. Scenario 0 represents the base case analyzed in Section \ref{subsec:results}.

\begin{table}[h!]
\ra{0.9}
\caption{Scenarios for the MA case}
\label{table: scenarios}
\scalebox{0.5}{
\begin{tabular}{llllll}
\toprule
\textbf{Scenario}   & \textbf{\# patients} & \textbf{Sync.} & \textbf{Patient types}             & \textbf{Period} & \textbf{Tot. ann. LFO}\\ \midrule \addlinespace[0.05cm]
\textbf{0} & 7000                        & 77\%                            & All patients                          & 4 months             & \euro \space  -130,874.47 
         \\
\textbf{1}                 & 7000                        & 77\%                            & All patients                          & 6 months   &  \euro \space -106,929.38 
                   \\
\textbf{2}                 & 7000                        & 77\%                            & Only patients with type $k_0$ & 4 months             &  \euro \space   -190,511.50 
        \\
\textbf{3}                 & 7000                        & 77\%                            & Only patients with type $k_0$  & 6 months         &  \euro \space     -179,020.58 
          \\
\textbf{4}                 & 7000                        & 87\%                            & All patients                          & 4 months     &  \euro \space     -124,559.17 
              \\
\textbf{5}                 & 7000                        & 87\%                            & All patients                          & 6 months      &  \euro \space   -103,014.38 
               \\
\textbf{6}                 & 7000                        & 87\%                            & Only patients with type $k_0$  & 4 months          &  \euro \space     -185,448.10 
         \\
\textbf{7}                 & 7000                        & 87\%                            & Only patients with type $k_0$  & 6 months           &  \euro \space     -131,990.02 
        \\
\textbf{8}                 & 7000                        & 100\%                           & All patients                          & 4 months    &  \euro \space   -117,366.55 
                 \\
\textbf{9}                 & 7000                        & 100\%                           & All patients                          & 6 months    &  \euro \space \space  -  97,687.64 
                 \\
\textbf{10}                & 7000                        & 100\%                           & Only patients with type $k_0$  & 4 months          &  \euro \space   -180,120.28 
           \\
\textbf{11}                & 7000                        & 100\%                           & Only patients with type $k_0$  & 6 months           &  \euro \space    -127,906.42 
         \\
\textbf{12}                & 10 000                      & 77\%                            & All patients                          & 4 months    &  \euro \space    -175,874.75 
                \\
\textbf{13}                & 10 000                      & 77\%                            & All patients                          & 6 months     &  \euro \space   -135,188.37 
                \\
\textbf{14}                & 10 000                      & 77\%                            & Only patients with type $k_0$  & 4 months        &  \euro \space   -206,118.67 
             \\
\textbf{15}                & 10 000                      & 77\%                            & Only patients with type $k_0$  & 6 months         &  \euro \space    -189,240.38 
           \\
\textbf{16}                & 10 000                      & 87\%                            & All patients                          & 4 months   &  \euro \space -139,635.38 
                 \\
\textbf{17}                & 10 000                      & 87\%                            & All patients                          & 6 months    &  \euro \space   -129,336.71 
                 \\
\textbf{18}                & 10 000                      & 87\%                            & Only patients with type $k_0$  & 4 months           &  \euro \space   -198,920.32 
          \\
\textbf{19}                & 10 000                      & 87\%                            & Only patients with type $k_0$  & 6 months           &  \euro \space   -184,463.56 
          \\
\textbf{20}                & 10 000                      & 100\%                           & All patients                          & 4 months    &  \euro \space\space   -80,187.73 
                 \\
\textbf{21}                & 10 000                      & 100\%                           & All patients                          & 6 months    &  \euro \space \space  -73,574.30 
                 \\
\textbf{22}                & 10 000                      & 100\%                           & Only patients with type $k_0$  & 4 months           &  \euro \space   -191,244.67 
          \\
\textbf{23}                & 10 000                      & 100\%                           & Only patients with type $k_0$  & 6 months           &  \euro \space   -179,309.32 
          \\ \bottomrule                        
\end{tabular}}
\end{table}

\begin{figure}[ht]
  \centering
  \subcaptionbox{KPI 1: Comparison of the annual total LFO for each scenario, including the LFO elements (transportation and handling costs and prescription line fee\label{fig: KPI1Scenarios}}{\includegraphics[width=0.45\textwidth]{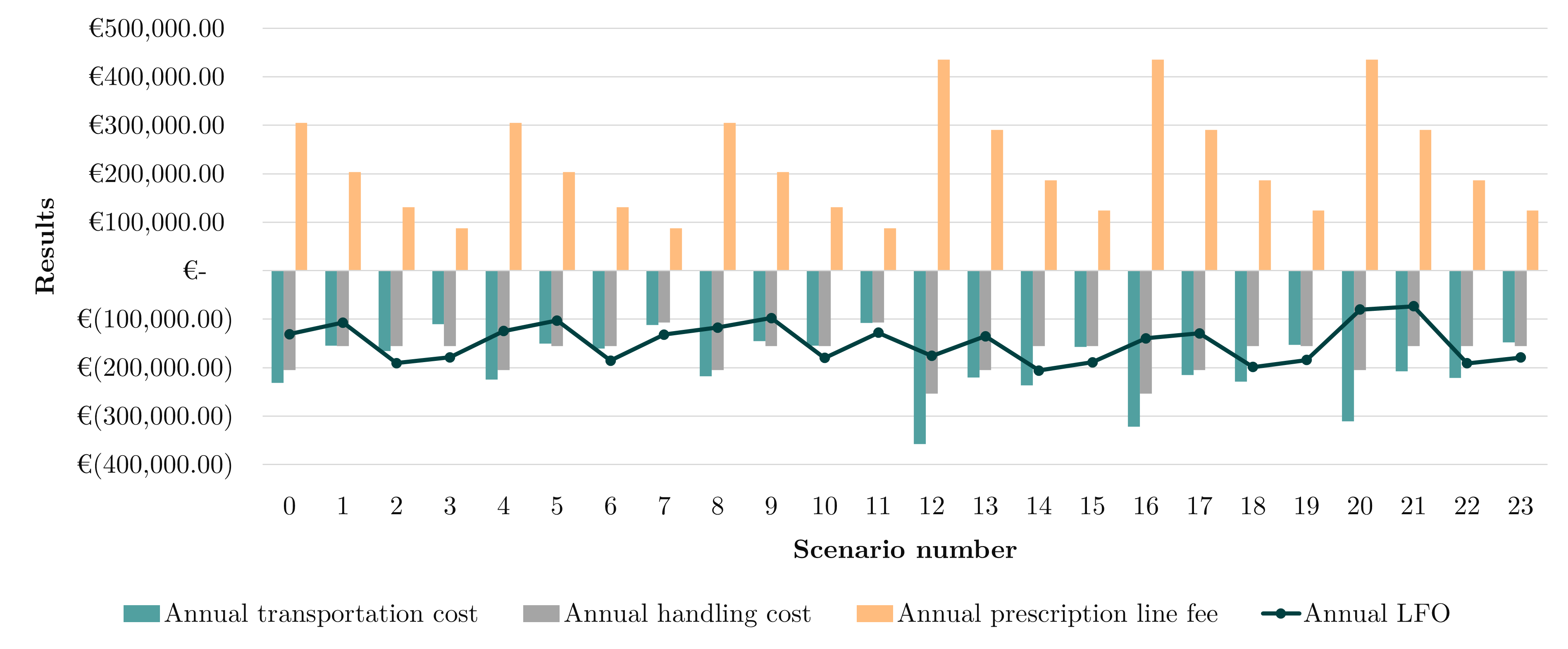}}
  \hfill
  \subcaptionbox{KPI 2: Comparison of the LFO per order batch, for each scenario, including the LFO elements (transportation and handling costs and prescription line fee\label{fig: KPI2Scenario}}{\includegraphics[width=0.45\textwidth]{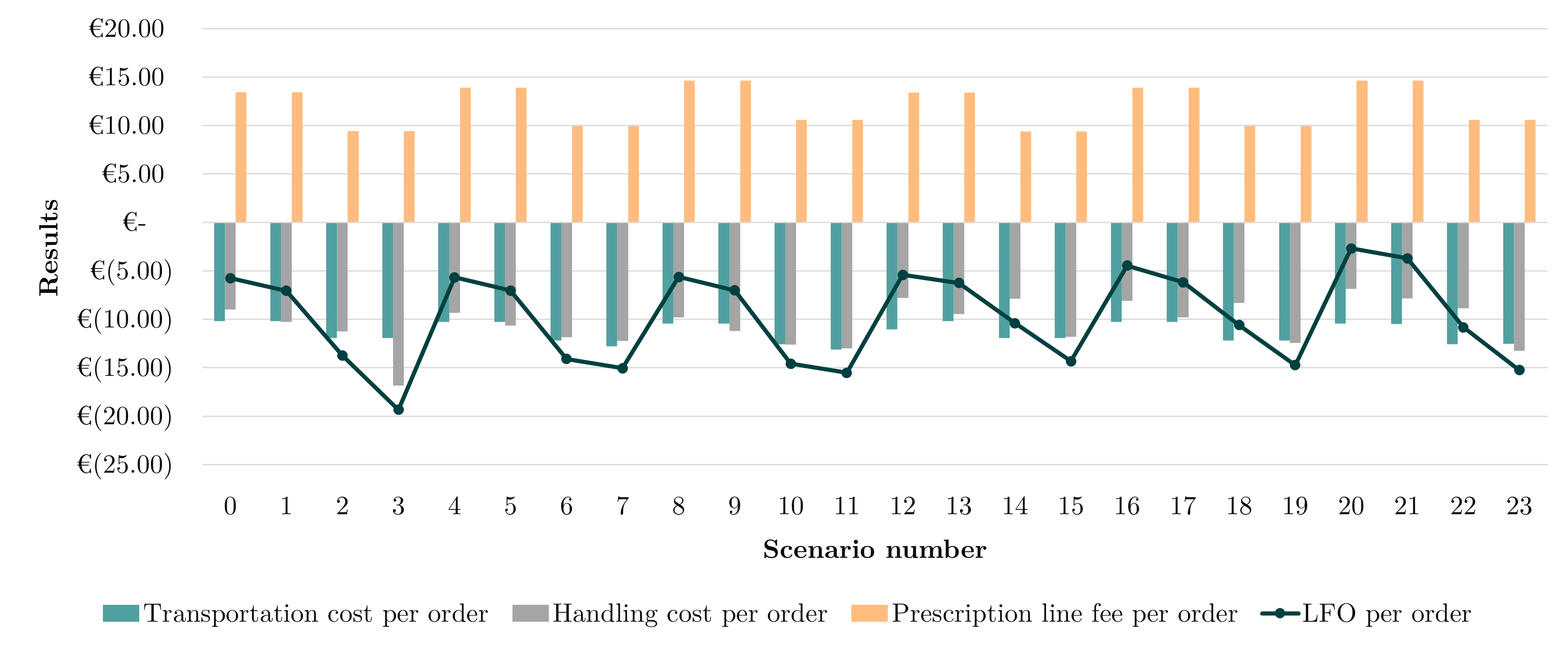}}
  \caption{Summary of the results for KPI 1 and KPI 2}
\end{figure}

To analyze the results of the scenarios (based on the KPIs defined in Section \ref{subsec:results}), we depict Figures \ref{fig: KPI1Scenarios} and \ref{fig: KPI2Scenario}. As displayed, the pharmacy experiences a net logistical loss in all scenarios. However, this outcome is not alarming, as the overall LFO of the pharmacy, considering all operations (beyond the scope of our model), is positive. The specific details of these operations and numbers are not presented due to privacy considerations. Among the scenarios, Scenario 21 stands out with the highest total annual LFO, while Scenario 20 exhibits the highest total annual LFO per order batch. An interesting pattern emerges when examining the differences among patient types. Figure \ref{fig: KPI2Scenario} demonstrates a consistently lower total annual LFO per order batch when only patients with medication type $k_0$ are served at the MA. Additionally, Figure \ref{fig: KPI1Scenarios} shows that scenarios with only patients of type $k_0$ yield a lower total annual LFO compared to their counterparts with mixed patient types.
The scenarios featuring only patients with medication type $k_0$ display a reduction in total costs, primarily due to a decrease in the number of patients. Most of these orders have no prescription line fee, so the associated costs are significantly lower. However, it is clear that exclusively serving patients with medication type $k_0$ is financially inefficient.
Therefore, the analysis concludes that limiting the MA's service to patients with medication type $k_0$ proves to be financially ineffective.



\begin{figure}[h!]
\centering
\resizebox{0.5\textwidth}{!}{
    \includegraphics[width = \textwidth]{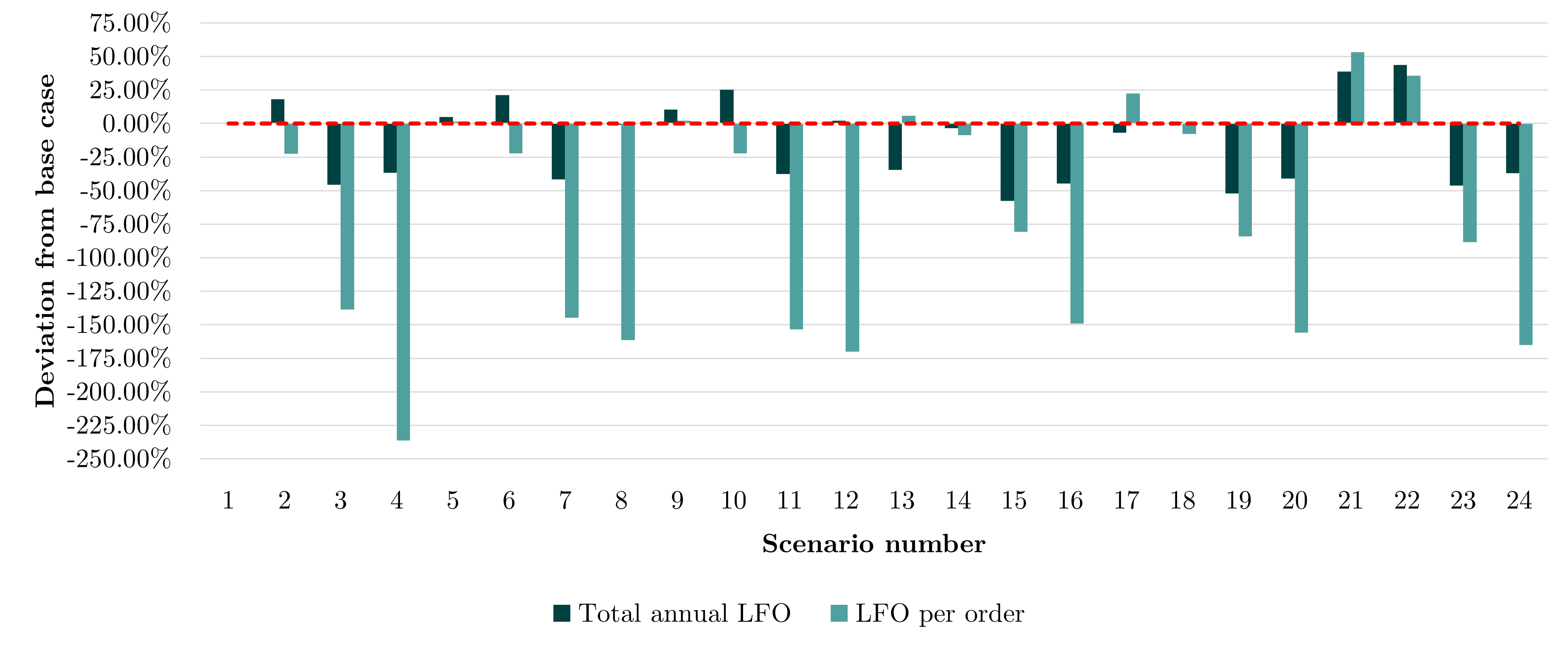}
    }
\caption{The influence on the base case result for each scenario: difference with KPIs 1 and 2 with respect to the base case (in percentage)}
\label{fig: ScenariosKPI12}
\end{figure}

Figure \ref{fig: ScenariosKPI12} provides an overview of the LFO differences compared to the base case, specifically focusing on the total annual LFO and the LFO per order batch (KPIs 1 and 2). Consistent with the findings from Figures \ref{fig: KPI1Scenarios} and \ref{fig: KPI2Scenario}, the patient type emerges as the most influential factor. Scenarios involving ``Only patients with type $k_0$" exhibit the largest decrease in both total annual LFO and LFO per order batch.
The number of patients has a variable impact on the LFO. The alternative number of patients (10,000 patients) demonstrates the most significant influence when combined with a different patient type scenario (scenario 14, yielding the lowest total annual LFO) or when combined with all alternative options (scenario 23, resulting in the lowest LFO per order batch for scenarios with 10,000 patients).
Altering the number of patients leads to a general decline in total annual LFO compared to the base case, although the LFO per order batch shows slight improvement. This outcome can be attributed to the fact that costs increase at a higher rate than total earnings, but the total annual LFO is divided among a larger number of orders.
Regarding the synchronization level, an increase in synchronization leads to higher total annual LFO and LFO per order batch. This correlation is logical since the prescription line fee remains constant across different synchronization levels (as indicated in Figure \ref{fig: KPI1Scenarios}), while costs rise with lower synchronization levels due to a greater number of orders.
Lastly, extending the order period from 4 to 6 months has a minor positive impact on the total annual LFO but adversely affects the LFO per order. This outcome can be attributed to the expectation of fewer orders when the order period is extended. Although total costs decrease, they are also spread across fewer orders.

\begin{figure}[ht]
  \centering
  \subcaptionbox{KPIs 3 and 4: Comparison of the number of order batches, divided into transportation cooling types\label{fig: ScenariosKPI34}}{\includegraphics[width=0.45\textwidth]{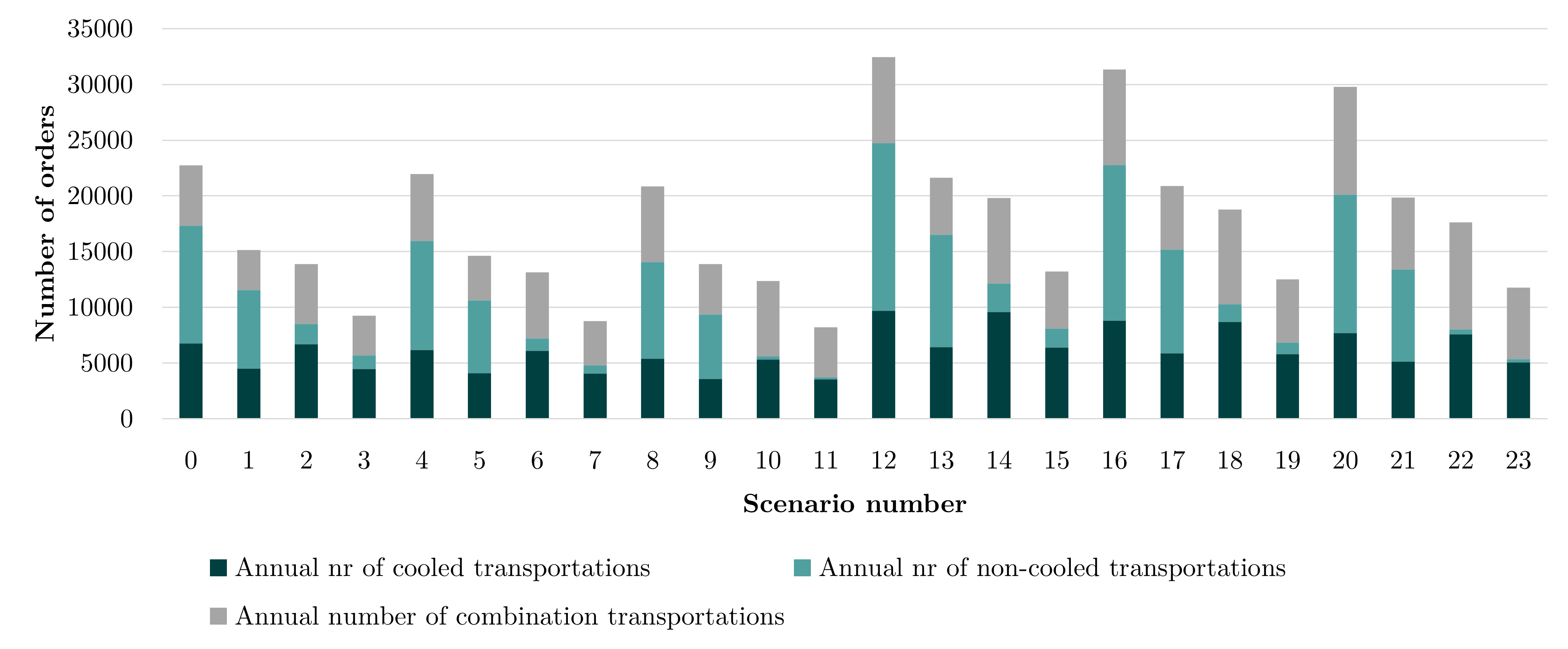}}
  \hfill
  \subcaptionbox{KPIs 3 and 5: Comparison of the number of order batches, divided into transportation types\label{fig: ScenariosKPI35}}{\includegraphics[width=0.45\textwidth]{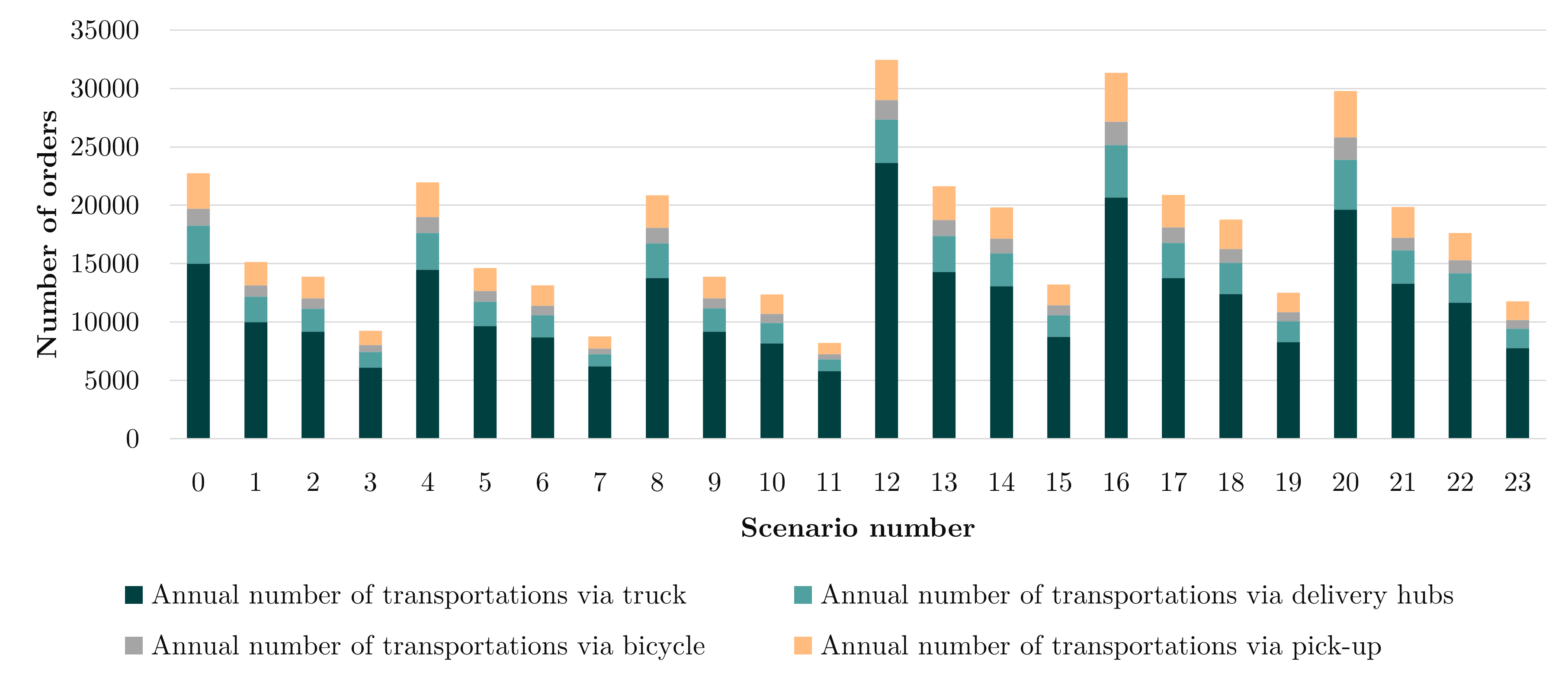}}
  \caption{Summary of the results for comparing KPIs}
\end{figure}

\begin{figure}[h!]
\centering
\resizebox{0.5\textwidth}{!}{

    \includegraphics[width = \textwidth]{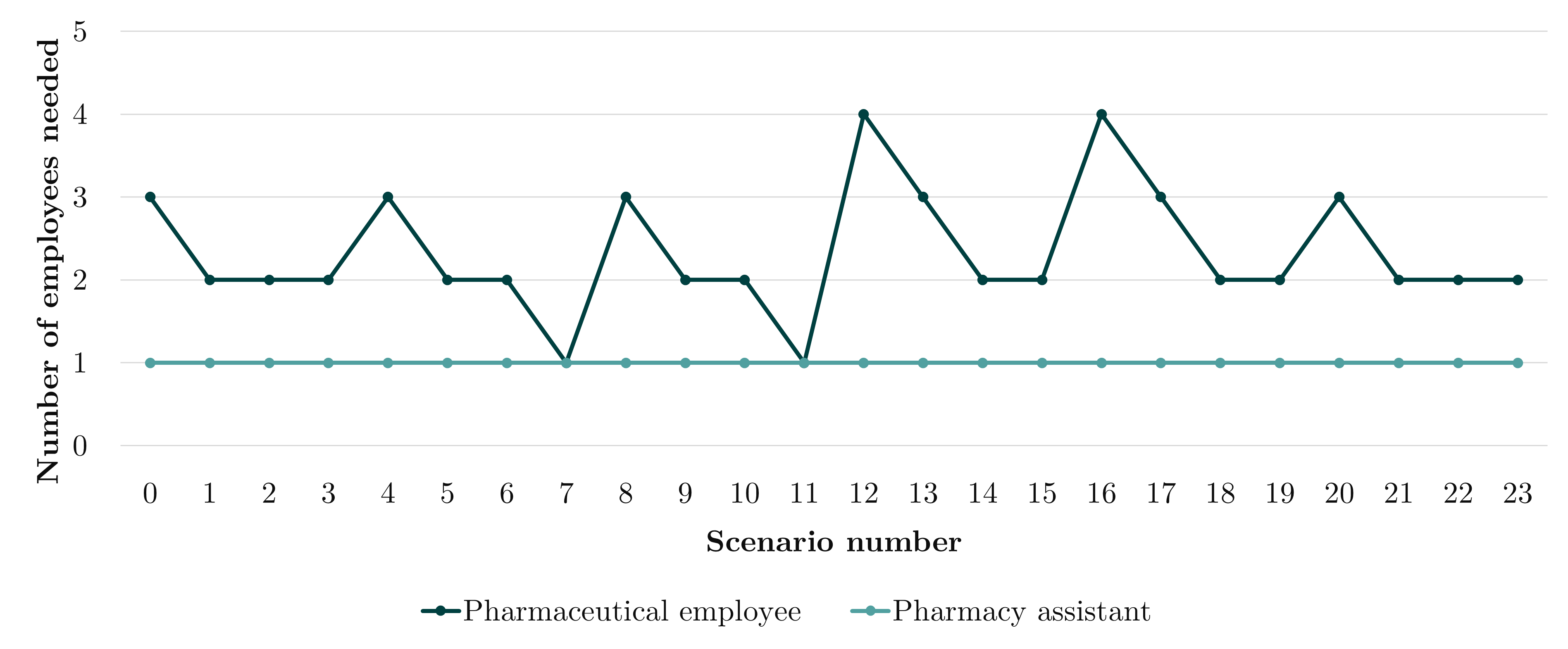}
    }
\caption{KPI 6: Comparison of the number of employees necessary in each scenario}
\label{fig: ScenariosKPI6}
\end{figure}

In the analysis of the decision KPIs, Figures \ref{fig: ScenariosKPI34} and \ref{fig: ScenariosKPI35} provide insights into the annual number of orders, which decreases under the following circumstances: (1) higher synchronization percentages, (2) exclusively serving patients with type $k_0$ medications, and (3) increasing the order period from 4 to 6 months. Conversely, the number of orders increases when the number of patients is increased.
Figure \ref{fig: ScenariosKPI34} illustrates that the ratio between cooling types differs primarily when patient types differ. This observation can be easily explained by the dominance of cooled medication types when exclusively serving patients with medication type $k_0$, as these medications require cooling during transportation.
Figure \ref{fig: ScenariosKPI35} indicates that the ratio between different delivery types remains relatively stable across scenarios. This suggests that the scenarios have limited influence on the investments in transportation modes (particularly cooling), mainly affecting the scale rather than the specific distribution of these investments.
Lastly, Figure \ref{fig: ScenariosKPI6} presents the number of employees required for the pharmacy's logistics. We observe a consistent pattern between the number of pharmaceutical employees and the number of order batches. However, the number of pharmacy technicians remains constant at one employee. This can be attributed to the fact that pharmacy technicians have limited involvement in the logistical processes of the MA, thereby not necessitating the addition of a second assistant in any scenario due to a low workload threshold.

The graphs provide valuable insights into the strategic decision-making process faced by the MA. From a financial perspective, scenarios 20 and 21 emerge as the most favorable choices. However, it is important to consider practical limitations when aiming for 100\% synchronization, as it may be challenging or even impossible to achieve in reality. Additionally, the MA's priorities may extend beyond maximizing LFO, encompassing factors such as total transportation cost or transportation costs per order.
In light of these considerations, it becomes crucial for management to carefully weigh these various factors and prioritize them according to the organization's specific goals and circumstances. The insights derived from this study can serve as a valuable resource to inform and support the decision-making process, helping management make well-informed strategic choices.







\subsection{Model variations} \label{subsec:modelvariations}
This section introduces two model variations aimed at enhancing the adaptability of the model to specific situations. These variations enable us to explore the influence of tactical-level decisions on the pharmacy's profitability. Furthermore, based on the MA case study, it becomes evident that the standard model may not always accurately reflect the reality of the logistics in the most realistic manner.

\subsubsection{Variation 1: relaxing the order numbers and order composition} 
In this model variation, we introduce changes to Equation \eqref{CUniqueMeds1}, resulting in Equations \eqref{SubEqVarUniqueMeds1} and \eqref{SubEqVarUniqueMeds2}. These modifications allow patients to either order their entire prescription every period $w$ or a portion of it. Equation \eqref{SubEqVarUniqueMeds1} ensures that the total demand over the entire time horizon is always met, ensuring that all medications are ordered at least once. On the other hand, Equation \eqref{SubEqVarUniqueMeds2} ensures that each order does not exceed the quantity prescribed. The objective of this model variation is to maximize the capacity utilization of order batches while adhering to the limitations imposed by the prescribed medications.
\begin{subequations}
    \begin{align}
        & \sum_{w\in W}o_{ckpw} \ge q_{ckp},   & \quad & \forall c\in C, k\in K, p\in P \label{SubEqVarUniqueMeds1}\\
        & o_{ckpw} \le  q_{ckp},  &\quad & \forall c\in C, k\in K, p\in P, w \in W \label{SubEqVarUniqueMeds2}
    \end{align}
\end{subequations}

This model variation yields intriguing outcomes. One notable observation is the substantial increase in the number of orders. This phenomenon can be attributed to the fact that, in the case of the MA, the prescription line fees tend to outweigh the associated order costs. However, it remains a pertinent question whether such frequent ordering is practical or viable for the pharmacy. Nevertheless, this insight provides valuable information for other types of pharmacies seeking to evaluate the trade-off between prescription line fees (or other revenue streams) and distribution costs.

\begin{figure}[ht]
  \centering
  \subcaptionbox{The annual \# orders of the base model vs. the variation model \label{fig: NrOrdersModelVariation}}{\includegraphics[width=0.45\textwidth]{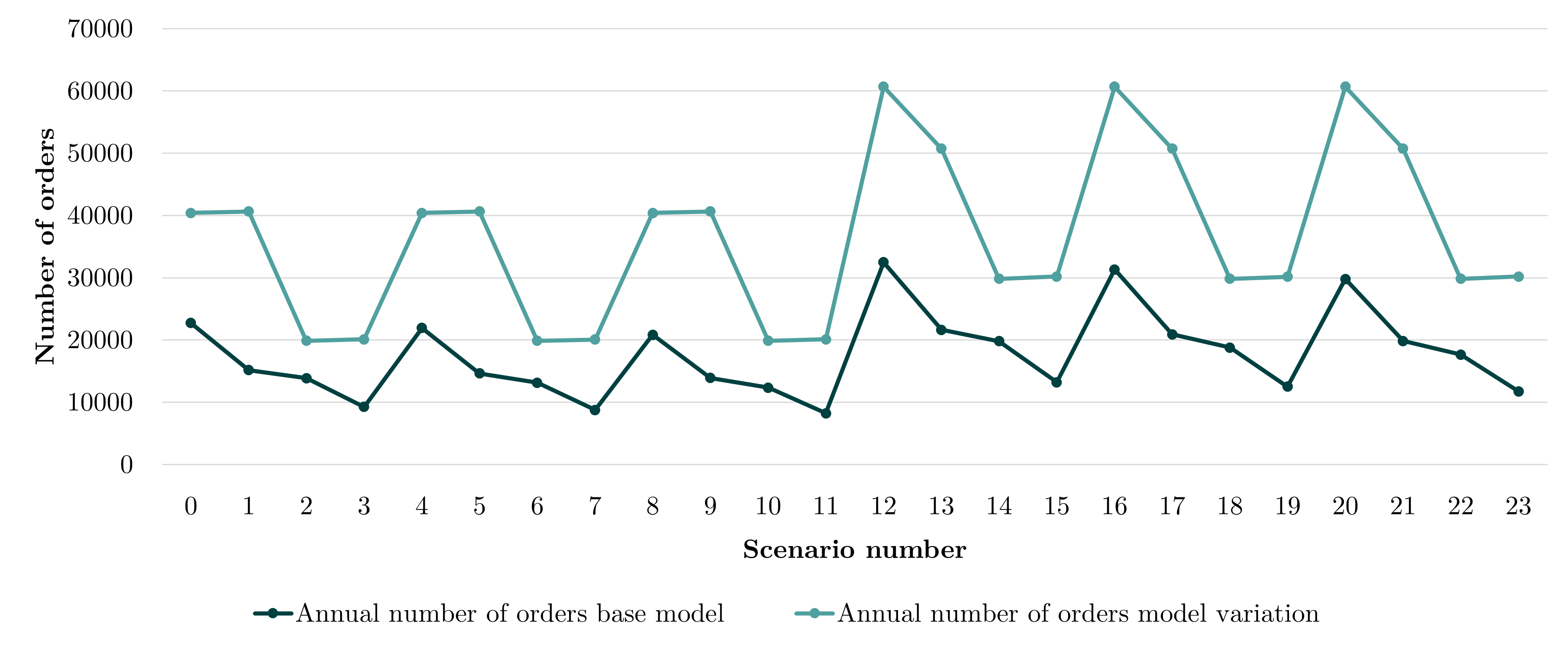}}
  \hfill
  \subcaptionbox{The annual LFO vs. the base model vs. the variation model \label{fig: ProfitModelVariation}}{\includegraphics[width=0.45\textwidth]{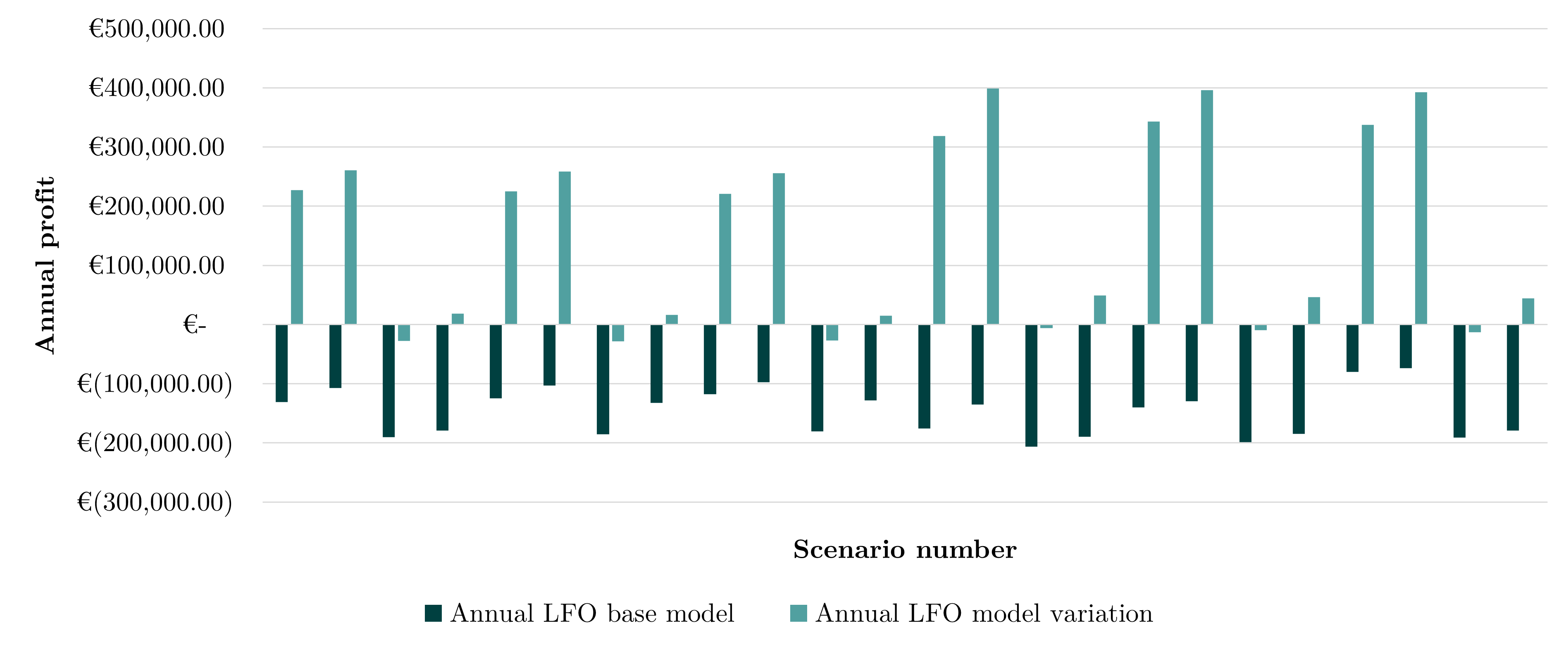}}
  \caption{Comparison between the base model and the variation}
\end{figure}

Figure \ref{fig: NrOrdersModelVariation} illustrates that, in most cases, the annual number of orders approximately doubles, accompanied by a more than threefold increase in LFO compared to the base model. This substantial growth can be attributed to the benefits of ordering more frequently in this model variation. The model intelligently evaluates whether the prescription line fee of a patient can offset the associated order costs, resulting in a positive LFO per order. 
To capitalize on these results, the pharmacy can adopt a strategy that encourages patients with multiple medications of type $k = k_1$ (which incur a non-zero prescription line fee) to order as frequently as possible, while simultaneously discouraging other patient types with a net loss from placing frequent orders.

\subsubsection{Variation 2: employee hours instead of number of employees}

This variation focuses on the number of employees required for logistical operations. According to KPI 6 and Figure \ref{fig: ScenariosKPI6}, we observe instances where up to four logistical employees are needed in the MA. The current model is designed to add an additional employee to the solution when the maximum number of hours for a specific employee type per period, denoted as $\theta_{ew}$, is exceeded. This approach is reasonable to prevent employees from consistently working overtime and to account for fluctuations in workload.

\begin{figure}[!ht]
\centering
\resizebox{0.45\textwidth}{!}{
    \includegraphics[width = \textwidth]{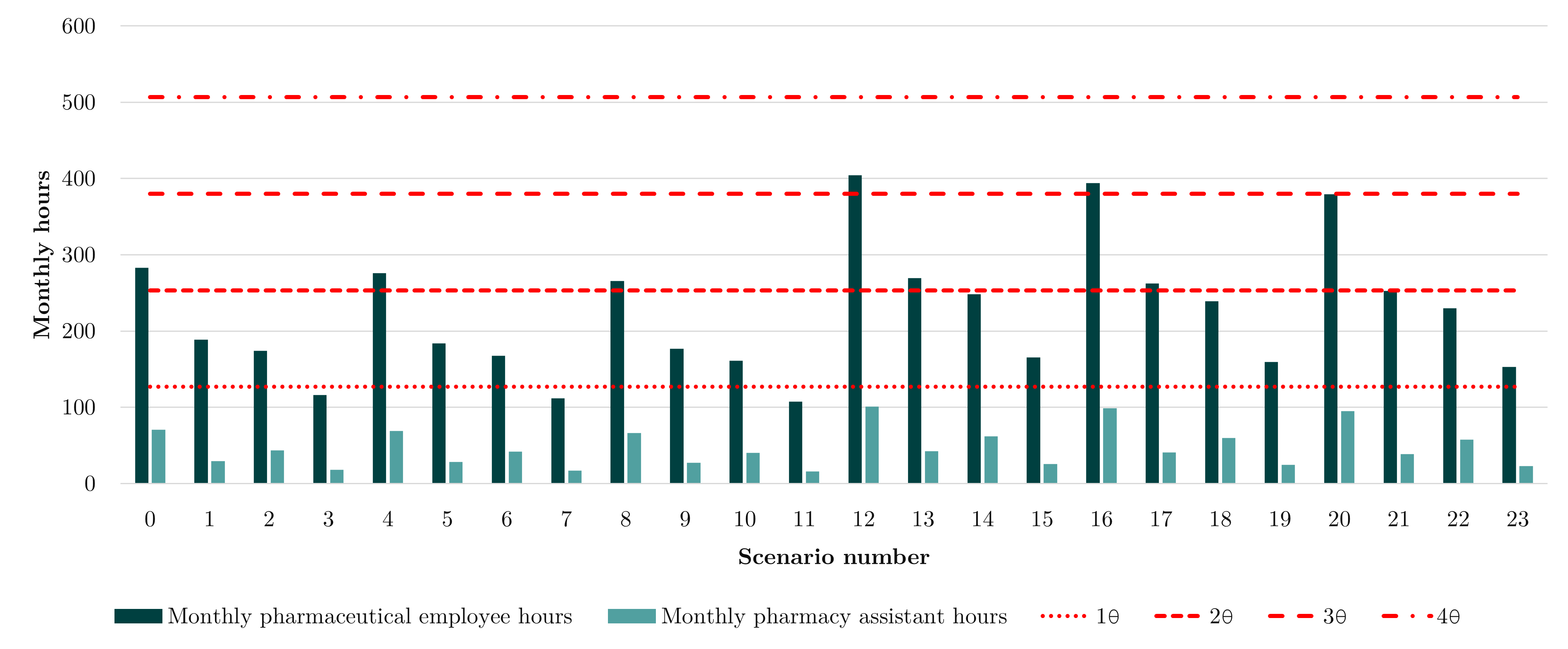}
    }
\caption{Number of hours needed of every employee type for each scenario}
\label{fig: NecessaryEmployeeHours}
\end{figure}

However, it is important to investigate the extent to which the $\theta_{ew}$ threshold is exceeded in each scenario. Figure \ref{fig: NecessaryEmployeeHours} displays the number of hours required for each employee type compared to the thresholds set for one to four employees (assuming the same value of $\theta_{ew}$ applies to all $e$ and $w$ in the MA case). We argue that the current strategy might not be the best one for the MA, as it results in a significant amount of idle time for employees overall. Moreover, in the MA case, employees are involved in tasks beyond logistics. They dedicate a portion of their hours to other aspects of pharmacy operations, such as attending to patients or working in the inpatient pharmacy. To address these concerns, we propose an alternative approach that eliminates the need for parameter $\theta_{ew}$ and instead focuses on calculating the necessary hours directly. This approach allows the pharmacy to have more flexibility in allocating these hours among their employees.

To implement this approach, we modify the definition of variable $m_{ew}$ to represent a non-integer variable, denoting ``the number of hours of employee $e$ needed in period $w$." Consequently, we remove variable $M_e$ and parameter $\theta_{ew}$ from the model. Additionally, we adjust the definition of parameter $s_e$ to represent ``the hourly salary of employee $e$."

By adopting this alternative formulation, the model provides the necessary hours required for each employee in each period, enabling the pharmacy to make informed decisions regarding the allocation of work hours among its employees. This approach offers greater flexibility and acknowledges that employees have responsibilities beyond logistics, accommodating their involvement in other pharmacy tasks such as patient care and inpatient pharmacy duties. 
The objective function is changed to Equation \eqref{ObjFuncVar2} to include the number of hours multiplied by the hourly wage of each employee $e$.

\begin{flalign}
    &\text{max } \sum_{c\in C} \sum_{k\in K}\sum_{p\in P}\sum_{w\in W} o_{ckpw} f_{k} \rho_{p}- \sum_{a\in A} \sum_{d\in D} \sum_{p\in P}\sum_{w\in W}  x_{adpw}t_{ad}\rho_{p}- \sum_{e\in E} \sum_{w\in W} m_{ew}s_{e}\label{ObjFuncVar2}
\end{flalign}
Constraint \eqref{CHoursPerWorker} is changed to Equation \eqref{CHoursPerWorkerVar2} and Constraint \eqref{CTotalWorkers} is removed from the model.
\begin{align}
 & \sum_{a\in A}\sum_{d\in D}\sum_{p\in P} x_{adpw}u_{ae}\rho_{p} \le m_{ew}    &\quad &  \forall w \in W \label{CHoursPerWorkerVar2}
\end{align}

This model variation provides valuable insights into scenarios where pharmacy employees are involved in tasks beyond logistics. By excluding hours dedicated to non-logistical responsibilities such as direct patient care or administrative work, the model accurately captures the workload specifically related to medication distribution. This approach enhances the accuracy of handling cost estimations and ensures the model's effectiveness in optimizing the logistical aspect of pharmacy operations.

It is important to note that this variation does not directly address the question of how many employees to hire, as it focuses on workload allocation rather than immediate staffing solutions. However, for pharmacies willing to make informed decisions about workload distribution among their employees, the output of this variation model can offer valuable insights into the necessary hours required for efficient logistical operations. By understanding the optimal allocation of labor hours, pharmacies can enhance operational efficiency and make informed decisions regarding resource utilization.

\subsection{Discussion} \label{subsec:methods statistics}
We model the MILP for the PCSCM-MSDM problem in Python 3.8 using the Python MIP package, and the Gurobi 10.0 solver. To evaluate the performance of the model, we conducted 200 runs of the base case (scenario 0) using the same data. This allowed us to perform a statistical analysis of the running times and optimality gap.

In all 200 runs, the Gurobi solver reported an optimality gap of no more than $1.8 * 10^{-13}\%$. The average running time for the base case was 3.99 seconds, with a 99\% confidence interval of [3.79, 4.19] seconds. These results demonstrate that the model runs efficiently and achieves near-optimality.

The base model consists of 14,410 variables. To assess the feasibility of the exact solution method using Gurobi for larger problem instances, we tested different set sizes. Table \ref{table:sidebyside} presents the number of indices and variables for three tested scenarios. The scenario numbers in the first column correspond to the first scenario (scenario 0) in Table \ref{table:sidebyside}, and the letters a, b, and c represent alternative variations in the number of variables.
Table \ref{table:sidebyside} provides information on the number of variables for the three scenarios. The variations in set $P$ are based on different classifications of the same patients, as described in Section \ref{subsec:data}. In scenario 0(a), the base case, there are 225 patient types. Scenario 0(b) includes 476 patient types without aggregating outlier patients, while scenario 0(c) involves separate classifications for each patient, resulting in 6783 patient types. The largest set, scenario 0(c), is 2912.64\% larger than the base scenario, with 434,122 variables compared to 14,410 variables in the base case, as shown in Table \ref{table:sidebyside}. These variations in set sizes allow us to examine the scalability and performance of the model for larger problem instances.

\begin{table}[h]
\captionsetup{singlelinecheck=false} 
\caption{Sizes of sets in different P size scenarios and Number of variables in different scenarios}
\label{table:sidebyside}
\begin{minipage}{0.2\linewidth}
\centering
\scalebox{0.7}{
\begin{tabular}{llllllll}
\toprule \addlinespace[0.1cm]
       & \textbf{Sets}  &            &            &            &            &            &            \\ \midrule \addlinespace[0.05cm]
\textbf{Scenario nr} & $\boldsymbol{A}$ & $\boldsymbol{C}$ & $\boldsymbol{D}$ & $\boldsymbol{E}$ & $\boldsymbol{K}$ & $\boldsymbol{P}$ & $\boldsymbol{W}$ \\ \midrule \addlinespace[0.1cm]
0(a)                      & 3          & 2          & 4          & 2          & 2          & 225        & 4          \\
0(b)                      & 3          & 2          & 4          & 2          & 2          & 476        & 4          \\
0(c)                      & 3          & 2          & 4          & 2          & 2          & 6,783       & 4          \\
\bottomrule
\end{tabular}}
\end{minipage}
\hfill
\begin{minipage}{0.6\linewidth}
\centering
\scalebox{0.7}{
\begin{tabular}{lllllll}
\toprule
\textbf{Scenario nr}    & $\boldsymbol{x_{adpw}}$  & $\boldsymbol{o_{ckpw}}$    & $\boldsymbol{m_{ew}}$ & $\boldsymbol{M_e}$   & \textbf{Total}  & \textbf{\% increase} \\	\midrule	\addlinespace[0.1cm]
0(a)                     & 10,800              & 3,600                 & 8                     & 2                 & 14,410             & 0\%                                \\
0(b)                     & 22,848              & 7,616                 & 8                     & 2                 & 30,474             & 111\%                              \\
0(c)                     & 325,584             & 108,528               & 8                     & 2                 & 434,122            & 2,913\%                             \\
\bottomrule
\end{tabular}}
\end{minipage}
\end{table}

Table \ref{tab: optimality gaps} shows the maximum optimality gap and running time 99\% confidence interval per scenario for $n=200$ runs. We see that case 0(b) has the largest maximum optimality gap of 0.0000753, which means that the solution is almost guaranteed to be near optimal using the Gurobi solver at all times. The largest running time found is 95.28 seconds for scenario 0(c), the scenario with the largest number of variables.

\begin{table}[H]
\caption{Maximum optimality gap and confidence interval of running time for different P size scenarios}
\label{tab: optimality gaps}
\scalebox{0.8}{
\begin{tabular}{lllll}
\toprule
\textbf{Scenario nr} & \textbf{Max. gap} & \textbf{Max. running time} &\textbf{99\% CI running time} & \textbf{LFO}\\ \midrule	\addlinespace[0.05cm]
0(a)               & $1.93*10^{-16}$              & 10.1 sec   & [3.79, 4.19]           &  \geneuro - 130,874.47     \\
0(b)               & $7.53*10^{-5}$               & 8.4 sec    & [3.50, 3.86]           &  \geneuro - 126,170.53    \\
0(c)               & $1.40*10^{-5}$               & 95.3 sec   & [76.83, 78.47]         &  \geneuro - 125,165,35 \\    \bottomrule
\end{tabular}}
\end{table}
We conclude that an alternative solution algorithm is not necessary at this time, since we see that increasing the model size with 2,912.64\% increases the maximum solution time to a very reasonable time still while maintaining a very small optimality gap. Additionally, the LFO stays approximately the same (within 5\% of the original solution of scenario 0(a)). Therefore, aggregating the classification of patients does not significantly sacrifice the quality of the solutions. We expect that for other pharmacies, the sizes of sets $K$, $E$, $P$, and $W$ have the most risk of being different than the set sizes in the MA case presented in this paper. If a pharmacy wishes to apply the model presented in this paper and has many medication types $k$, employee types $e$, patient types $p$, and periods $w$, such that the number of variables is significantly larger than the number of variables of scenario 0(c) (434,122), the model has a risk of being too slow, especially when the pharmacy wishes to run multiple scenarios. In this case, alternative solution methods may be considered, albeit this would be an extreme case that we do not expect to occur.

\section{Conclusion}\label{sec: Conclusion}

In this study, we addressed the Pharmaceutical Cold Supply Chain Management Considering Medication Synchronization and Different Delivery Modes (PCSCM-MDSM) problem. Our goal was to examine the effects of delivery methods, medication synchronization, and cooling requirements on distribution costs and revenues. To tackle this problem, we developed a Mixed-Integer Linear Programming (MILP) model. This model aims to optimize financial results (LFO) by determining the optimal number and composition of orders, as well as the appropriate number of employees while considering patient demand, cooling restrictions, and employee constraints.
To evaluate the effectiveness of our model, we utilized real-world data from the Maartensapotheek (MA), an outpatient pharmacy affiliated with a Dutch hospital. We implemented and solved the MILP model using the Gurobi 10.0 solver within the Python-MIP package. Our solver successfully provided near-optimal solutions for all instances of the MA case (an optimality gap of $1.83 * 10^{-15}$ for the base scenario), demonstrating the efficacy of our approach.
Furthermore, we conducted tests with two variations of the model and explored two alternative model sizes. The maximum observed optimality gap across these experiments was $7.53 * 10^{-5}$. Computation times of alternative model sizes are never more than 100 seconds, which ensures the model runs within a very reasonable time. These findings demonstrate the robustness and reliability of our model in addressing the PCSCM-MDSM problem and its ability to generate high-quality solutions in a reasonable time.

Despite the promising outcomes of this research, there exist certain limitations that warrant acknowledgment and present avenues for future investigation. Firstly, the current model overlooks dosage and package sizes. While package size seldom poses challenges for couriers in practice, variations in delivery modes may introduce packaging size issues. For instance, when patients order less frequently, such as every six months, the packages may become larger if the number of orders per period remains constant. Further research should explore the feasibility of potentially increasing order sizes for both the pharmacy and the chosen delivery methods. Additionally, the model assumes an absence of stock shortages in the pharmacy, reflecting the infrequency of such shortages in the MA. However, this may not hold true for different pharmacies, where shortages could occur. In such cases, pharmacies often opt for more frequent patient orders to facilitate smaller batches. This adjustment would significantly influence the model's results, as it would generate a larger number of orders for certain medications. Consequently, investigating the impact of shortages on the model and incorporating appropriate adaptations becomes imperative for future research. Moreover, the model's current formulation encompasses two types of cooling restrictions and three types of transportation cooling restrictions, one for each medication cooling restriction, as well as a combination type. Should the model be employed in a pharmacy with additional cooling types, the constraints outlined in Section \ref{subsec: mathematicalformulation} would require modification. Another potential area for future research, as mentioned in Section \ref{subsec:methods statistics}, lies in exploring alternative solution methods that may prove advantageous for larger problem instances (it is not useful for our case as we could solve the problem to optimality in a matter of seconds). Lastly, although not necessary for our case, expanding the number of scenarios tested in future research would provide a broader evaluation of the model's performance across diverse environments, enhancing its practical applicability.

We also note that the current model is of deterministic nature, as the variability of input data is deemed to be negligible by the problem owner. However, the problem can be imposed on stochastic elements corresponding to the input data. For instance, factors such as demand arrival patterns, frequency, and the composition of order batches may exhibit variability over time, influenced by different patient types. To address this dynamic environment more comprehensively, future research should consider incorporating stochastic models and exploring appropriate solution methods, such as simulation modeling. By adopting such approaches, decision-makers can obtain valuable insights and make informed decisions in the face of uncertainty.


\nocite{*}

\bibliographystyle{itor}
\bibliography{references1}

\section*{Appendix}
This section contains the expanded versions of values reported in Table \ref{tab: summaryqrhosigma}.
\small 
\renewcommand{\arraystretch}{0.8} 

\begin{longtable}[h]{lllllllllll}
\caption{Complete values of $q_{ckp}$: the number of unique medicines with cooling type $c$, medication type $k$, and for patient type $p$}
\label{tab:completeq} \\
\hline
& \multicolumn{4}{l}{\textbf{Nr unique meds $\boldsymbol{q_{ckp}}$}} & & & \multicolumn{4}{l}{\textbf{Nr unique meds $\boldsymbol{q_{ckp}}$}} \\
\cline{2-6} \cline{7-11}
\textbf{Type $\boldsymbol{p}$} & $c_0 ; k_0$       & $c_0 ; k_1$  & $c_1 ; k_0$  & $c_1 ; k_1$ & & \textbf{Type $\boldsymbol{p}$} & $c_0 ; k_0$       & $c_0 ; k_1$  & $c_1 ; k_0$  & $c_1 ; k_1$ \\
\cline{2-6} \cline{7-11}
\endfirsthead
\hline

\cline{1-5} \cline{7-11}
\multicolumn{11}{r}{{Continued on next page}} \\
\hline
\endfoot

\hline
\multicolumn{11}{r}{{Continued on next page}} \\
\cline{1-5} \cline{7-11}
\endfoot

\cline{1-5} \cline{7-11}
\multicolumn{11}{r}{{Continued on next page}} \\
\cline{1-5} \cline{7-11}
\endfoot

\cline{1-5} \cline{7-11}
\textbf{Type $\boldsymbol{p}$} & $c_0 ; k_0$       & $c_0 ; k_1$  & $c_1 ; k_0$  & $c_1 ; k_1$ & & \textbf{Type $\boldsymbol{p}$} & $c_0 ; k_0$       & $c_0 ; k_1$  & $c_1 ; k_0$  & $c_1 ; k_1$ \\
\cline{1-5} \cline{7-11}
\endhead

\cline{1-5} \cline{7-11}
\multicolumn{11}{r}{{Continued on next page}} \\
\hline
\endfoot

\cline{1-5} \cline{7-11}
\multicolumn{11}{r}{{End of table}} \\
\hline
\endlastfoot
        
$p_{0}$	&	1	&	0	&	0	&	0	&	&
$p_{1}$	&	1	&	0	&	0	&	0		\\
$p_{2}$	&	0	&	0	&	0	&	1	&	&
$p_{3}$	&	0	&	0	&	0	&	1		\\
$p_{4}$	&	0	&	0	&	0	&	2	&	&
$p_{5}$	&	1	&	0	&	0	&	1		\\
$p_{6}$	&	0	&	0	&	0	&	2	&	&
$p_{7}$	&	0	&	0	&	0	&	3		\\
$p_{8}$	&	0	&	0	&	0	&	2	&	&
$p_{9}$	&	1	&	0	&	0	&	1		\\
$p_{10}$	&	1	&	0	&	0	&	2	&	&
$p_{11}$	&	0	&	0	&	0	&	3		\\
$p_{12}$	&	0	&	0	&	0	&	3	&	&
$p_{13}$	&	1	&	0	&	0	&	2		\\
$p_{14}$	&	0	&	0	&	0	&	4	&	&
$p_{15}$	&	0	&	0	&	0	&	4		\\
$p_{16}$	&	2	&	0	&	0	&	0	&	&
$p_{17}$	&	1	&	0	&	0	&	3		\\
$p_{18}$	&	0	&	0	&	0	&	4	&	&
$p_{19}$	&	1	&	0	&	0	&	1		\\
$p_{20}$	&	1	&	0	&	0	&	1	&	&
$p_{21}$	&	2	&	0	&	0	&	0		\\
$p_{22}$	&	1	&	0	&	0	&	3	&	&
$p_{23}$	&	1	&	0	&	0	&	2		\\
$p_{24}$	&	2	&	0	&	0	&	0	&	&
$p_{25}$	&	0	&	0	&	0	&	2		\\
$p_{26}$	&	1	&	0	&	0	&	3	&	&
$p_{27}$	&	1	&	0	&	0	&	2		\\
$p_{28}$	&	0	&	0	&	0	&	5	&	&
$p_{29}$	&	1	&	0	&	0	&	3		\\
$p_{30}$	&	0	&	0	&	0	&	5	&	&
$p_{31}$	&	1	&	0	&	0	&	4		\\
$p_{32}$	&	1	&	0	&	0	&	4	&	&
$p_{33}$	&	0	&	0	&	0	&	5		\\
$p_{34}$	&	1	&	0	&	0	&	1	&	&
$p_{35}$	&	0	&	0	&	0	&	3		\\
$p_{36}$	&	1	&	0	&	0	&	2	&	&
$p_{37}$	&	1	&	0	&	0	&	4		\\
$p_{38}$	&	1	&	0	&	0	&	2	&	&
$p_{39}$	&	0	&	0	&	0	&	3		\\
$p_{40}$	&	0	&	0	&	0	&	6	&	&
$p_{41}$	&	0	&	0	&	0	&	4		\\
$p_{42}$	&	1	&	0	&	0	&	2	&	&
$p_{43}$	&	0	&	0	&	1	&	0		\\
$p_{44}$	&	0	&	0	&	0	&	5	&	&
$p_{45}$	&	2	&	0	&	0	&	1		\\
$p_{46}$	&	1	&	0	&	0	&	5	&	&
$p_{47}$	&	2	&	0	&	0	&	1		\\
$p_{48}$	&	1	&	0	&	0	&	3	&	&
$p_{49}$	&	0	&	1	&	0	&	0		\\
$p_{50}$	&	1	&	0	&	0	&	5	&	&
$p_{51}$	&	0	&	0	&	0	&	6		\\
$p_{52}$	&	1	&	0	&	0	&	3	&	&
$p_{53}$	&	0	&	0	&	0	&	4		\\
$p_{54}$	&	1	&	0	&	0	&	5	&	&
$p_{55}$	&	1	&	0	&	0	&	1		\\
$p_{56}$	&	0	&	0	&	0	&	6	&	&
$p_{57}$	&	2	&	0	&	0	&	1		\\
$p_{58}$	&	1	&	0	&	0	&	3	&	&
$p_{59}$	&	1	&	0	&	0	&	5		\\
$p_{60}$	&	0	&	0	&	0	&	7	&	&
$p_{61}$	&	0	&	0	&	0	&	3		\\
$p_{62}$	&	1	&	1	&	0	&	0	&	&
$p_{63}$	&	0	&	0	&	0	&	6		\\
$p_{64}$	&	0	&	0	&	0	&	7	&	&
$p_{65}$	&	2	&	0	&	0	&	2		\\
$p_{66}$	&	2	&	0	&	0	&	3	&	&
$p_{67}$	&	1	&	0	&	1	&	0		\\
$p_{68}$	&	1	&	0	&	0	&	5	&	&
$p_{69}$	&	1	&	0	&	0	&	4		\\
$p_{70}$	&	1	&	0	&	1	&	1	&	&
$p_{71}$	&	2	&	0	&	0	&	2		\\
$p_{72}$	&	1	&	0	&	0	&	4	&	&
$p_{73}$	&	1	&	0	&	0	&	3		\\
$p_{74}$	&	1	&	0	&	0	&	2	&	&
$p_{75}$	&	0	&	0	&	0	&	4		\\
$p_{76}$	&	1	&	0	&	0	&	6	&	&
$p_{77}$	&	2	&	0	&	0	&	1		\\
$p_{78}$	&	2	&	0	&	0	&	3	&	&
$p_{79}$	&	1	&	1	&	0	&	0		\\
$p_{80}$	&	2	&	0	&	0	&	4	&	&
$p_{81}$	&	1	&	0	&	0	&	6		\\
$p_{82}$	&	2	&	0	&	0	&	0	&	&
$p_{83}$	&	0	&	0	&	1	&	0		\\
$p_{84}$	&	2	&	0	&	0	&	2	&	&
$p_{85}$	&	2	&	0	&	0	&	1		\\
$p_{86}$	&	1	&	0	&	0	&	1	&	&
$p_{87}$	&	2	&	0	&	0	&	0		\\
$p_{88}$	&	0	&	0	&	1	&	1	&	&
$p_{89}$	&	1	&	0	&	0	&	6		\\
$p_{90}$	&	0	&	0	&	0	&	5	&	&
$p_{91}$	&	0	&	0	&	0	&	2		\\
$p_{92}$	&	1	&	0	&	0	&	4	&	&
$p_{93}$	&	0	&	0	&	0	&	4		\\
$p_{94}$	&	1	&	0	&	1	&	0	&	&
$p_{95}$	&	0	&	0	&	0	&	7		\\
$p_{96}$	&	2	&	0	&	0	&	3	&	&
$p_{97}$	&	2	&	0	&	0	&	1		\\
$p_{98}$	&	1	&	0	&	0	&	2	&	&
$p_{99}$	&	0	&	0	&	0	&	5		\\
$p_{100}$	&	2	&	0	&	0	&	1	&	&
$p_{101}$	&	1	&	0	&	0	&	6		\\
$p_{102}$	&	1	&	0	&	0	&	7	&	&
$p_{103}$	&	2	&	0	&	0	&	2		\\
$p_{104}$	&	1	&	0	&	1	&	1	&	&
$p_{105}$	&	1	&	0	&	0	&	3		\\
$p_{106}$	&	0	&	0	&	0	&	3	&	&
$p_{107}$	&	2	&	0	&	0	&	1		\\
$p_{108}$	&	1	&	0	&	0	&	7	&	&
$p_{109}$	&	0	&	0	&	0	&	7		\\
$p_{110}$	&	0	&	0	&	1	&	2	&	&
$p_{111}$	&	2	&	0	&	0	&	4		\\
$p_{112}$	&	0	&	0	&	0	&	9	&	&
$p_{113}$	&	2	&	0	&	0	&	5		\\
$p_{114}$	&	2	&	0	&	0	&	5	&	&
$p_{115}$	&	2	&	0	&	0	&	3		\\
$p_{116}$	&	2	&	0	&	0	&	2	&	&
$p_{117}$	&	2	&	0	&	0	&	4		\\
$p_{118}$	&	0	&	0	&	0	&	8	&	&
$p_{119}$	&	1	&	0	&	0	&	4		\\
$p_{120}$	&	1	&	0	&	1	&	0	&	&
$p_{121}$	&	1	&	0	&	0	&	2		\\
$p_{122}$	&	0	&	0	&	1	&	2	&	&
$p_{123}$	&	1	&	0	&	0	&	1		\\
$p_{124}$	&	1	&	0	&	0	&	5	&	&
$p_{125}$	&	0	&	0	&	0	&	8		\\
$p_{126}$	&	2	&	0	&	0	&	3	&	&
$p_{127}$	&	0	&	0	&	0	&	5		\\
$p_{128}$	&	2	&	0	&	0	&	4	&	&
$p_{129}$	&	2	&	0	&	0	&	0		\\
$p_{130}$	&	0	&	0	&	1	&	1	&	&
$p_{131}$	&	0	&	0	&	1	&	3		\\
$p_{132}$	&	2	&	0	&	0	&	0	&	&
$p_{133}$	&	1	&	0	&	0	&	2		\\
$p_{134}$	&	1	&	0	&	0	&	5	&	&
$p_{135}$	&	1	&	0	&	1	&	2		\\
$p_{136}$	&	2	&	0	&	0	&	6	&	&
$p_{137}$	&	0	&	0	&	0	&	9		\\
$p_{138}$	&	1	&	0	&	0	&	2	&	&
$p_{139}$	&	1	&	0	&	0	&	4		\\
$p_{140}$	&	1	&	0	&	0	&	4	&	&
$p_{141}$	&	1	&	0	&	0	&	6		\\
$p_{142}$	&	0	&	1	&	0	&	1	&	&
$p_{143}$	&	1	&	0	&	1	&	0		\\
$p_{144}$	&	1	&	0	&	0	&	4	&	&
$p_{145}$	&	1	&	0	&	1	&	0		\\
$p_{146}$	&	1	&	0	&	0	&	3	&	&
$p_{147}$	&	2	&	0	&	0	&	4		\\
$p_{148}$	&	0	&	0	&	0	&	2	&	&
$p_{149}$	&	0	&	1	&	0	&	5		\\
$p_{150}$	&	0	&	0	&	0	&	2	&	&
$p_{151}$	&	1	&	0	&	0	&	1		\\
$p_{152}$	&	1	&	0	&	1	&	3	&	&
$p_{153}$	&	1	&	1	&	0	&	3		\\
$p_{154}$	&	1	&	0	&	0	&	1	&	&
$p_{155}$	&	1	&	0	&	0	&	7		\\
$p_{156}$	&	1	&	0	&	1	&	1	&	&
$p_{157}$	&	0	&	0	&	1	&	3		\\
$p_{158}$	&	2	&	0	&	0	&	2	&	&
$p_{159}$	&	1	&	0	&	1	&	4		\\
$p_{160}$	&	0	&	0	&	1	&	2	&	&
$p_{161}$	&	2	&	0	&	0	&	2		\\
$p_{162}$	&	0	&	0	&	0	&	8	&	&
$p_{163}$	&	0	&	0	&	1	&	3		\\
$p_{164}$	&	1	&	0	&	1	&	2	&	&
$p_{165}$	&	1	&	0	&	0	&	3		\\
$p_{166}$	&	2	&	0	&	0	&	2	&	&
$p_{167}$	&	0	&	1	&	0	&	1		\\
$p_{168}$	&	1	&	0	&	0	&	4	&	&
$p_{169}$	&	0	&	1	&	0	&	1		\\
$p_{170}$	&	1	&	0	&	0	&	5	&	&
$p_{171}$	&	1	&	0	&	0	&	2		\\
$p_{172}$	&	0	&	0	&	0	&	2	&	&
$p_{173}$	&	2	&	0	&	0	&	5		\\
$p_{174}$	&	1	&	0	&	0	&	2	&	&
$p_{175}$	&	1	&	0	&	1	&	2		\\
$p_{176}$	&	0	&	0	&	1	&	5	&	&
$p_{177}$	&	1	&	0	&	0	&	5		\\
$p_{178}$	&	1	&	0	&	0	&	8	&	&
$p_{179}$	&	1	&	1	&	0	&	0		\\
$p_{180}$	&	1	&	1	&	0	&	0	&	&
$p_{181}$	&	1	&	0	&	1	&	3		\\
$p_{182}$	&	2	&	0	&	0	&	3	&	&
$p_{183}$	&	1	&	0	&	0	&	7		\\
$p_{184}$	&	1	&	0	&	0	&	5	&	&
$p_{185}$	&	2	&	0	&	0	&	5		\\
$p_{186}$	&	0	&	0	&	1	&	4	&	&
$p_{187}$	&	1	&	0	&	0	&	6		\\
$p_{188}$	&	1	&	0	&	0	&	9	&	&
$p_{189}$	&	1	&	0	&	0	&	7		\\
$p_{190}$	&	0	&	0	&	1	&	1	&	&
$p_{191}$	&	0	&	0	&	0	&	8		\\
$p_{192}$	&	1	&	0	&	1	&	2	&	&
$p_{193}$	&	1	&	0	&	0	&	14		\\
$p_{194}$	&	2	&	0	&	0	&	3	&	&
$p_{195}$	&	1	&	0	&	0	&	8		\\
$p_{196}$	&	1	&	0	&	0	&	9	&	&
$p_{197}$	&	0	&	0	&	1	&	1		\\
$p_{198}$	&	2	&	0	&	0	&	2	&	&
$p_{199}$	&	1	&	1	&	0	&	1		\\
$p_{200}$	&	2	&	0	&	0	&	3	&	&
$p_{201}$	&	1	&	0	&	0	&	7		\\
$p_{202}$	&	0	&	0	&	1	&	3	&	&
$p_{203}$	&	0	&	0	&	0	&	10		\\
$p_{204}$	&	1	&	0	&	1	&	2	&	&
$p_{205}$	&	2	&	0	&	0	&	1		\\
$p_{206}$	&	1	&	1	&	0	&	0	&	&
$p_{207}$	&	2	&	0	&	0	&	6		\\
$p_{208}$	&	1	&	0	&	1	&	0	&	&
$p_{209}$	&	1	&	0	&	0	&	8		\\
$p_{210}$	&	0	&	0	&	0	&	9	&	&
$p_{211}$	&	1	&	1	&	0	&	2		\\
$p_{212}$	&	2	&	0	&	0	&	4	&	&
$p_{213}$	&	0	&	0	&	1	&	4		\\
$p_{214}$	&	1	&	1	&	1	&	5	&	&
$p_{215}$	&	0	&	1	&	0	&	2		\\
$p_{216}$	&	0	&	2	&	0	&	0	&	&
$p_{217}$	&	0	&	0	&	1	&	5		\\
$p_{218}$	&	1	&	0	&	1	&	1	&	&
$p_{219}$	&	0	&	1	&	0	&	4		\\
$p_{220}$	&	0	&	1	&	0	&	3	&	&
$p_{221}$	&	1	&	0	&	1	&	3		\\
$p_{222}$	&	1	&	1	&	0	&	1	&	&
$p_{223}$	&	1	&	0	&	1	&	5		\\
$p_{224}$	&	0	&	1	&	0	&	6	&	&

\end{longtable}

\normalsize
\small 
\renewcommand{\arraystretch}{0.95} 

\begin{table}[H]
\caption{Complete values of $\rho_{p}$: the number of patients of type $p$ and $\sigma_{p}$: the minimum number of orders per time horizon of patient type $p$}
    \label{tab: completesigma}
    \scalebox{0.77}{
\begin{tabular}{lll|lll|lll|lll|lll|lll}
\hline
\textbf{$\boldsymbol{p}$} & \textbf{\# $\boldsymbol{\rho_p}$} & \textbf{\# $\boldsymbol{\sigma_p}$} &
\textbf{$\boldsymbol{p}$} & \textbf{\# $\boldsymbol{\rho_p}$} & \textbf{\# $\boldsymbol{\sigma_p}$} &
\textbf{$\boldsymbol{p}$} & \textbf{\# $\boldsymbol{\rho_p}$} & \textbf{\# $\boldsymbol{\sigma_p}$} &
\textbf{$\boldsymbol{p}$} & \textbf{\# $\boldsymbol{\rho_p}$} & \textbf{\# $\boldsymbol{\sigma_p}$} &
\textbf{$\boldsymbol{p}$} & \textbf{\# $\boldsymbol{\rho_p}$} & \textbf{\# $\boldsymbol{\sigma_p}$} &
\textbf{$\boldsymbol{p}$} & \textbf{\# $\boldsymbol{\rho_p}$} & \textbf{\# $\boldsymbol{\sigma_p}$} \\
\hline
$\boldsymbol{p_{0}}$ & 935 & 1 & $\boldsymbol{p_{40}}$ & 29 & 1 & $\boldsymbol{p_{80}}$ & 11 & 1 & $\boldsymbol{p_{120}}$ & 6 & 1 & $\boldsymbol{p_{160}}$ & 4 & 1 & $\boldsymbol{p_{200}}$ & 3 & 1 \\
$\boldsymbol{p_{1}}$ & 553 & 1 & $\boldsymbol{p_{41}}$ & 29 & 1 & $\boldsymbol{p_{81}}$ & 11 & 1 & $\boldsymbol{p_{121}}$ & 6 & 3 & $\boldsymbol{p_{161}}$ & 4 & 1 & $\boldsymbol{p_{201}}$ & 3 & 1 \\
$\boldsymbol{p_{2}}$ & 539 & 1 & $\boldsymbol{p_{42}}$ & 28 & 2 & $\boldsymbol{p_{82}}$ & 10 & 1 & $\boldsymbol{p_{122}}$ & 6 & 1 & $\boldsymbol{p_{162}}$ & 4 & 1 & $\boldsymbol{p_{202}}$ & 3 & 1 \\
$\boldsymbol{p_{3}}$ & 318 & 1 & $\boldsymbol{p_{43}}$ & 31 & 1 & $\boldsymbol{p_{83}}$ & 11 & 1 & $\boldsymbol{p_{123}}$ & 6 & 1 & $\boldsymbol{p_{163}}$ & 3 & 2 & $\boldsymbol{p_{203}}$ & 2 & 1 \\
$\boldsymbol{p_{4}}$ & 261 & 1 & $\boldsymbol{p_{44}}$ & 25 & 1 & $\boldsymbol{p_{84}}$ & 10 & 1 & $\boldsymbol{p_{124}}$ & 6 & 2 & $\boldsymbol{p_{164}}$ & 3 & 2 & $\boldsymbol{p_{204}}$ & 2 & 1 \\
$\boldsymbol{p_{5}}$ & 259 & 1 & $\boldsymbol{p_{45}}$ & 25 & 1 & $\boldsymbol{p_{85}}$ & 10 & 1 & $\boldsymbol{p_{125}}$ & 6 & 1 & $\boldsymbol{p_{165}}$ & 3 & 3 & $\boldsymbol{p_{205}}$ & 2 & 1 \\
$\boldsymbol{p_{6}}$ & 208 & 1 & $\boldsymbol{p_{46}}$ & 24 & 1 & $\boldsymbol{p_{86}}$ & 10 & 1 & $\boldsymbol{p_{126}}$ & 6 & 1 & $\boldsymbol{p_{166}}$ & 3 & 3 & $\boldsymbol{p_{206}}$ & 2 & 1 \\
$\boldsymbol{p_{7}}$ & 180 & 1 & $\boldsymbol{p_{47}}$ & 22 & 1 & $\boldsymbol{p_{87}}$ & 10 & 1 & $\boldsymbol{p_{127}}$ & 6 & 2 & $\boldsymbol{p_{167}}$ & 3 & 1 & $\boldsymbol{p_{207}}$ & 2 & 1 \\
$\boldsymbol{p_{8}}$ & 183 & 1 & $\boldsymbol{p_{48}}$ & 22 & 2 & $\boldsymbol{p_{88}}$ & 9 & 1 & $\boldsymbol{p_{128}}$ & 6 & 1 & $\boldsymbol{p_{168}}$ & 3 & 3 & $\boldsymbol{p_{208}}$ & 2 & 1 \\
$\boldsymbol{p_{9}}$ & 157 & 1 & $\boldsymbol{p_{49}}$ & 22 & 1 & $\boldsymbol{p_{89}}$ & 9 & 1 & $\boldsymbol{p_{129}}$ & 6 & 1 & $\boldsymbol{p_{169}}$ & 3 & 1 & $\boldsymbol{p_{209}}$ & 1 & 1 \\
$\boldsymbol{p_{10}}$ & 154 & 1 & $\boldsymbol{p_{50}}$ & 21 & 1 & $\boldsymbol{p_{90}}$ & 9 & 2 & $\boldsymbol{p_{130}}$ & 6 & 1 & $\boldsymbol{p_{170}}$ & 3 & 3 & $\boldsymbol{p_{210}}$ & 1 & 1 \\
$\boldsymbol{p_{11}}$	&	147	&	1	&	$\boldsymbol{p_{51}}$	&	20	&	1	&	$\boldsymbol{p_{91}}$	&	9	&	1	&	$\boldsymbol{p_{131}}$	&	5	&	1	&	$\boldsymbol{p_{171}}$	&	3	&	3	&	$\boldsymbol{p_{211}}$	&	1	&	1	\\
$\boldsymbol{p_{12}}$	&	124	&	1	&	$\boldsymbol{p_{52}}$	&	18	&	2	&	$\boldsymbol{p_{92}}$	&	9	&	2	&	$\boldsymbol{p_{132}}$	&	5	&	1	&	$\boldsymbol{p_{172}}$	&	3	&	1	&	$\boldsymbol{p_{212}}$	&	1	&	1	\\
$\boldsymbol{p_{13}}$	&	110	&	1	&	$\boldsymbol{p_{53}}$	&	18	&	2	&	$\boldsymbol{p_{93}}$	&	9	&	2	&	$\boldsymbol{p_{133}}$	&	5	&	2	&	$\boldsymbol{p_{173}}$	&	3	&	2	&	$\boldsymbol{p_{213}}$	&	1	&	1	\\
$\boldsymbol{p_{14}}$	&	93	&	1	&	$\boldsymbol{p_{54}}$	&	18	&	1	&	$\boldsymbol{p_{94}}$	&	9	&	1	&	$\boldsymbol{p_{134}}$	&	5	&	2	&	$\boldsymbol{p_{174}}$	&	3	&	3	&	$\boldsymbol{p_{214}}$	&	1	&	1	\\
$\boldsymbol{p_{15}}$	&	97	&	1	&	$\boldsymbol{p_{55}}$	&	18	&	1	&	$\boldsymbol{p_{95}}$	&	9	&	1	&	$\boldsymbol{p_{135}}$	&	5	&	2	&	$\boldsymbol{p_{175}}$	&	3	&	1	&	$\boldsymbol{p_{215}}$	&	1	&	1	\\
$\boldsymbol{p_{16}}$	&	93	&	1	&	$\boldsymbol{p_{56}}$	&	18	&	1	&	$\boldsymbol{p_{96}}$	&	8	&	1	&	$\boldsymbol{p_{136}}$	&	5	&	2	&	$\boldsymbol{p_{176}}$	&	3	&	1	&	$\boldsymbol{p_{216}}$	&	1	&	1	\\
$\boldsymbol{p_{17}}$	&	92	&	1	&	$\boldsymbol{p_{57}}$	&	17	&	2	&	$\boldsymbol{p_{97}}$	&	8	&	2	&	$\boldsymbol{p_{137}}$	&	5	&	1	&	$\boldsymbol{p_{177}}$	&	3	&	2	&	$\boldsymbol{p_{217}}$	&	1	&	1	\\
$\boldsymbol{p_{18}}$	&	77	&	1	&	$\boldsymbol{p_{58}}$	&	17	&	2	&	$\boldsymbol{p_{98}}$	&	8	&	3	&	$\boldsymbol{p_{138}}$	&	5	&	2	&	$\boldsymbol{p_{178}}$	&	3	&	1	&	$\boldsymbol{p_{218}}$	&	1	&	1	\\
$\boldsymbol{p_{19}}$	&	73	&	1	&	$\boldsymbol{p_{59}}$	&	16	&	2	&	$\boldsymbol{p_{99}}$	&	8	&	2	&	$\boldsymbol{p_{139}}$	&	5	&	2	&	$\boldsymbol{p_{179}}$	&	3	&	1	&	$\boldsymbol{p_{219}}$	&	1	&	1	\\
$\boldsymbol{p_{20}}$	&	71	&	1	&	$\boldsymbol{p_{60}}$	&	16	&	1	&	$\boldsymbol{p_{100}}$	&	8	&	3	&	$\boldsymbol{p_{140}}$	&	5	&	2	&	$\boldsymbol{p_{180}}$	&	3	&	1	&	$\boldsymbol{p_{220}}$	&	1	&	1	\\
$\boldsymbol{p_{21}}$	&	65	&	1	&	$\boldsymbol{p_{61}}$	&	15	&	2	&	$\boldsymbol{p_{101}}$	&	8	&	2	&	$\boldsymbol{p_{141}}$	&	5	&	1	&	$\boldsymbol{p_{181}}$	&	3	&	1	&	$\boldsymbol{p_{221}}$	&	1	&	1	\\
$\boldsymbol{p_{22}}$	&	62	&	1	&	$\boldsymbol{p_{62}}$	&	15	&	1	&	$\boldsymbol{p_{102}}$	&	8	&	1	&	$\boldsymbol{p_{142}}$	&	5	&	1	&	$\boldsymbol{p_{182}}$	&	3	&	2	&	$\boldsymbol{p_{222}}$	&	1	&	1	\\
$\boldsymbol{p_{23}}$	&	61	&	1	&	$\boldsymbol{p_{63}}$	&	15	&	1	&	$\boldsymbol{p_{103}}$	&	8	&	2	&	$\boldsymbol{p_{143}}$	&	4	&	1	&	$\boldsymbol{p_{183}}$	&	3	&	2	&	$\boldsymbol{p_{223}}$	&	1	&	1	\\
$\boldsymbol{p_{24}}$	&	60	&	1	&	$\boldsymbol{p_{64}}$	&	15	&	1	&	$\boldsymbol{p_{104}}$	&	7	&	2	&	$\boldsymbol{p_{144}}$	&	4	&	3	&	$\boldsymbol{p_{184}}$	&	3	&	2	&	$\boldsymbol{p_{224}}$	&	1	&	1	\\
$\boldsymbol{p_{25}}$	&	54	&	1	&	$\boldsymbol{p_{65}}$	&	14	&	1	&	$\boldsymbol{p_{105}}$	&	7	&	2	&	$\boldsymbol{p_{145}}$	&	4	&	1	&	$\boldsymbol{p_{185}}$	&	3	&	1	&	&	\\			
$\boldsymbol{p_{26}}$	&	53	&	1	&	$\boldsymbol{p_{66}}$	&	14	&	2	&	$\boldsymbol{p_{106}}$	&	7	&	2	&	$\boldsymbol{p_{146}}$	&	4	&	3	&	$\boldsymbol{p_{186}}$	&	3	&	1	&	&	\\			
$\boldsymbol{p_{27}}$	&	45	&	1	&	$\boldsymbol{p_{67}}$	&	14	&	1	&	$\boldsymbol{p_{107}}$	&	7	&	2	&	$\boldsymbol{p_{147}}$	&	4	&	2	&	$\boldsymbol{p_{187}}$	&	3	&	2	&	&	\\				
$\boldsymbol{p_{28}}$	&	45	&	1	&	$\boldsymbol{p_{68}}$	&	14	&	1	&	$\boldsymbol{p_{108}}$	&	7	&	1	&	$\boldsymbol{p_{148}}$	&	4	&	1	&	$\boldsymbol{p_{188}}$	&	3	&	1	&	&	\\				
$\boldsymbol{p_{29}}$	&	41	&	1	&	$\boldsymbol{p_{69}}$	&	13	&	2	&	$\boldsymbol{p_{109}}$	&	7	&	1	&	$\boldsymbol{p_{149}}$	&	4	&	1	&	$\boldsymbol{p_{189}}$	&	3	&	2	&	&	\\				
$\boldsymbol{p_{30}}$	&	41	&	1	&	$\boldsymbol{p_{70}}$	&	13	&	1	&	$\boldsymbol{p_{110}}$	&	7	&	1	&	$\boldsymbol{p_{150}}$	&	4	&	1	&	$\boldsymbol{p_{190}}$	&	3	&	1	&	&	\\				
$\boldsymbol{p_{31}}$	&	40	&	1	&	$\boldsymbol{p_{71}}$	&	13	&	1	&	$\boldsymbol{p_{111}}$	&	6	&	1	&	$\boldsymbol{p_{151}}$	&	4	&	1	&	$\boldsymbol{p_{191}}$	&	3	&	1	&	&	\\				
$\boldsymbol{p_{32}}$	&	40	&	1	&	$\boldsymbol{p_{72}}$	&	13	&	1	&	$\boldsymbol{p_{112}}$	&	6	&	1	&	$\boldsymbol{p_{152}}$	&	4	&	2	&	$\boldsymbol{p_{192}}$	&	3	&	1	&	&	\\				
$\boldsymbol{p_{33}}$	&	40	&	1	&	$\boldsymbol{p_{73}}$	&	11	&	2	&	$\boldsymbol{p_{113}}$	&	6	&	1	&	$\boldsymbol{p_{153}}$	&	4	&	1	&	$\boldsymbol{p_{193}}$	&	3	&	2	&	&	\\			
$\boldsymbol{p_{34}}$	&	39	&	1	&	$\boldsymbol{p_{74}}$	&	11	&	2	&	$\boldsymbol{p_{114}}$	&	6	&	1	&	$\boldsymbol{p_{154}}$	&	4	&	1	&	$\boldsymbol{p_{194}}$	&	3	&	2	&	&	\\				
$\boldsymbol{p_{35}}$	&	37	&	2	&	$\boldsymbol{p_{75}}$	&	11	&	2	&	$\boldsymbol{p_{115}}$	&	6	&	2	&	$\boldsymbol{p_{155}}$	&	4	&	1	&	$\boldsymbol{p_{195}}$	&	3	&	1	&	&	\\				
$\boldsymbol{p_{36}}$	&	36	&	2	&	$\boldsymbol{p_{76}}$	&	11	&	1	&	$\boldsymbol{p_{116}}$	&	6	&	2	&	$\boldsymbol{p_{156}}$	&	4	&	2	&	$\boldsymbol{p_{196}}$	&	3	&	1	&	&	\\				
$\boldsymbol{p_{37}}$	&	34	&	1	&	$\boldsymbol{p_{77}}$	&	11	&	2	&	$\boldsymbol{p_{117}}$	&	6	&	2	&	$\boldsymbol{p_{157}}$	&	4	&	1	&	$\boldsymbol{p_{197}}$	&	3	&	1	&	&	\\				
$\boldsymbol{p_{38}}$	&	29	&	2	&	$\boldsymbol{p_{78}}$	&	11	&	1	&	$\boldsymbol{p_{118}}$	&	6	&	1	&	$\boldsymbol{p_{158}}$	&	4	&	2	&	$\boldsymbol{p_{198}}$	&	3	&	2	&	&	\\				
$\boldsymbol{p_{39}}$	&	29	&	1	&	$\boldsymbol{p_{79}}$	&	11	&	1	&	$\boldsymbol{p_{119}}$	&	6	&	2	&	$\boldsymbol{p_{159}}$	&	4	&	1	&	$\boldsymbol{p_{199}}$	&	3	&	1	&	&	\\	\hline	\\		

\end{tabular}}
\end{table}

\end{document}